\documentclass[fontsize=10pt,paper=a4,headings= small,
    headsepline=true,abstract=true,twoside,]{scrartcl}

\KOMAoption{abstract}{true}
\setkomafont{section}{\large}
\setkomafont{subsection}{\normalsize}
\setkomafont{subsubsection}{\normalsize}
\setkomafont{paragraph}{\normalsize}
\setkomafont{caption}{\small}
\setkomafont{captionlabel}{\bfseries\usekomafont{caption}}

\usepackage[automark]{scrlayer-scrpage}
  \clearpairofpagestyles
  \ihead[]{\headmark}
\usepackage[inner=1.8cm, outer=3cm, top=3cm, bottom=4cm]{geometry}

\usepackage[T1]{fontenc}
\usepackage[utf8]{inputenc}
\usepackage{lmodern}
\usepackage[british]{babel}

\usepackage{graphicx}
\usepackage[export]{adjustbox}
\usepackage[dvipsnames]{xcolor}
\definecolor{darkred}{RGB}{139,0,0}
\definecolor{mediumblue}{RGB}{0,0,205}
\definecolor{forestgreen}{RGB}{34,139,34}

\usepackage[auth-sc-lg, noblocks]{authblk}

\usepackage{csquotes}
\usepackage[
    sorting=nty,
    sortcites=true,
    style=numeric-comp,
    doi=true,
    url=true,
    maxbibnames=4,
    backend=biber,
    giveninits=true,
    isbn=false,
    ]{biblatex}
\bibliography{literature}

\AtEveryBibitem{\clearfield{month}}
\AtEveryBibitem{\clearfield{day}}

\usepackage{amssymb, amsmath, amsthm}
\usepackage{mathtools}
\usepackage{commath}
\usepackage{nicefrac}
\usepackage{array}
\usepackage{multirow}

\usepackage{bm}
\usepackage{calrsfs}
\usepackage{stmaryrd}
\DeclareMathAlphabet{\pazocal}{OMS}{zplm}{m}{n}
\DeclareMathAlphabet{\pazocalbf}{OMS}{cmsy}{b}{n}
\usepackage{etoolbox}
\preto\subequations{\ifhmode\unskip\fi}

\usepackage[colorlinks=true]{hyperref}
\hypersetup{linkcolor=darkred, urlcolor=forestgreen, citecolor=mediumblue}
\usepackage{cleveref}

\usepackage{enumitem}
  \setlist[enumerate]{nosep, topsep=0pt, wide = 1em, leftmargin=*}
  \setlist[itemize]{nosep, topsep=3pt, wide = 1em, leftmargin=*}
\usepackage{booktabs}
\AtBeginEnvironment{tabular}{\small}

\newtheorem{form}{Formulation}
\newtheorem{Algorithm}{Algorithm}
\theoremstyle{remark}
\newtheorem{remark}{Remark}

\def\ep{\varepsilon}

\let\phivar\phi
\def\phi{\varphi}

\newcommand{\hcalA}{\hat{\cal A}}

\newcommand{\RR}{\mathbb{R}}
\newcommand{\PP}{\mathbb{P}}

\newcommand{\xb}{\bm{x}}
\newcommand{\xbref}{\hat{\xb}}
\newcommand{\nb}{\bm{n}}

\newcommand{\CR}{\pazocal{C}}

\renewcommand{\O}{\Omega}
\newcommand{\Oref}{\hat{\O}}

\newcommand{\FL}{\pazocal{F}}
\newcommand{\FLref}{\hat{\pazocal{F}}}
\newcommand{\SO}{\pazocal{S}}
\newcommand{\SOref}{\hat{\pazocal{S}}}
\newcommand{\IN}{\pazocal{I}}
\newcommand{\INref}{\hat{\pazocal{I}}}

\newcommand{\phiref}{\hat{\bm{\phivar}}}
\newcommand{\xiref}{\hat{\xi}}
\newcommand{\psiref}{\hat{\bm{\psi}}}

\newcommand{\fb}{\bm{f}}
\newcommand{\fbref}{\hat{\fb}}
\newcommand{\gb}{\bm{g}}

\newcommand{\ub}{\bm{u}}
\newcommand{\ubref}{\hat{\ub}}

\newcommand{\vb}{\bm{v}}
\newcommand{\vbref}{\hat{\vb}}

\newcommand{\wb}{\bm{w}}

\newcommand{\sigb}{\bm{\sigma}}
\newcommand{\sigbref}{\hat{\sigb}}
\newcommand{\eb}{\bm{e}}

\newcommand{\Vb}{\bm{V}}
\newcommand{\Vbref}{\hat{\Vb}}

\DeclareMathOperator{\COD}{cod}
\DeclareMathOperator{\TCV}{tcv}
\DeclareMathOperator{\meas}{meas}

\DeclareMathOperator{\diver}{div}
\newcommand{\nrm}[1]{\Vert#1\Vert}
\newcommand{\restr}[2]{{\left.\kern-\nulldelimiterspace#1\right|_{#2}}}

\title{\Large
A high-accuracy framework for phase-field fracture interface 
reconstructions with application to Stokes fluid-filled fracture
surrounded by an elastic medium}
\author[1,*]{Henry von Wahl} \author[2,**]{Thomas Wick} 

\affil[1]{Fakultät für Mathematik, Universität Wien, Oskar-Morgenstern-Platz, 1090 Wien, Austria.}
\affil[*]{E-Mail: \texttt{henry.vonwahl@ovgu.de}}
\affil[2]{Institut für Angewandte Mathematik, Leibniz Universität Hannover, Welfengarten 1, 30167 Hannover, Germany.}
\affil[**]{E-Mail: \texttt{thomas.wick@ifam.uni-hannover.de}}

\date{\small \today}

\begin{document}
\maketitle

\begin{abstract}

This work considers a Stokes flow in a deformable fracture interacting with a linear elastic medium. To this end, we employ a phase-field model to approximate the crack dynamics. Phase-field methods belong to interface-capturing approaches in which the interface is only given by a smeared zone. For multi-domain problems, the accuracy of the coupling conditions is, 
however, of utmost importance. Here, interface-tracking methods are preferred,
since the interface is resolved on mesh edges up to discretization errors, but it does not depend on the length scale parameter of some smeared zone. The key objective of this work is to construct a robust framework that computes first a crack path via the phase-field method (interface-capturing) and then does an interface-tracking reconstruction. We then discuss several approaches to reconstruct the Eulerian description of the open crack domain. This includes unfitted approaches where a level-set of the crack interface is constructed and an approach where the geometry is re-meshed. Using this reconstructed domain, we can compute the fluid-structure interaction problem between the fluid in the crack and the interacting solid. With the explicit mesh reconstruction of the two domains, we can then use an interface-tracking Arbitrary-Lagrangian-Eulerian (ALE) discretisation approach for the resulting fluid-structure interaction (FSI) problem. Our algorithmic procedure is realised in one final numerical
algorithm and one implementation. We substantiate our approach using several numerical examples based on Sneddon's benchmark and corresponding extensions to Stokes fluid-filled regimes.
\end{abstract}

\section{Introduction}
This work is devoted to a coupled multi-domain multi-physics problem between a fluid-filled, pressurised crack and the surrounding elastic solid in which the crack develops. Applications for this can be found in porous media problems such as natural and induced fractures, groundwater flow, nuclear waste management, and fluid-filled biomaterials.

A well-known approach for fluids, multi-physics and fracture propagation is the \emph{phase-field} approach. The so-called Stefan problem was subject in \cite{Plotnikov1992}.
A numerical analysis of phase-field in flow problems was established in \cite{FePro04},
for two-phase incompressible flows in \cite{ShenYang10}, and for multi-component 
flows in \cite{kim_2012}. Fully Eulerian 
phase-field methods for modelling fluid-structure interaction were investigated 
in \cite{Sun20141,MAO2023111903}. Phase-field methods for modelling 
tumour growth include for example \cite{GARCKE2021103192,FRITZ2021103331}.
Variational phase-field fracture, the most relevant for this paper, has been studied 
by various groups, as in \cite{BourFraMar00,KuMue10,BoVeScoHuLa12,AmGeraLoren15,HEIDER2018116,Wilson2016264,AlNoWiWr21_CAMWA,WheWiLee20,DiLiWiTy22}.
Finally, textbooks, monographs, and extended papers of phase-field methods in 
material sciences and fracture are for instance \cite{ProEl10,BourFraMar08,WNN20,Wic20}.

The most attractive features of this method are a fixed background mesh and the possibility for 
interfaces to move up to large deformations and topology changes, i.e., contact. A less favourable property in all phase-field methods 
is that the crack interface is smeared. 
While in many applications, this is no point of concern, 
this becomes a challenge when interface conditions need to 
be described \cite{Nguyen2015,GARCKE2021103192,MWW19}.
In some situations, this issue can be circumvented using mathematical formulae such as the Gauss divergence theorem. However, this elegant approach requires care due to the different boundaries (see \cite[Section 6.2]{MWW19}). Furthermore, such approaches are limited to `simple' physics on the interface. The extension of the phase-field approach to multi-physics with varying solutions and coefficients on the interface is challenging, and numerical solutions are highly dependent on the interface thickness (phase-field regularisation parameter, also known as length scale) and its relation to the spatial discretisation. 
Several ideas to reconstruct or approximate the interface have been proposed to deal with the problem of the smeared interface in a phase-field approach in the context of multi-physics problems, for example, by using additional displacement fields \cite{VeBo13} or explicit level-set functions \cite{Nguyen2015}. The latter idea was first applied to fluid-filled 
fractures in \cite{LWW17}, yielding a reliable method. However, this requires the computation 
of an additional Laplace-extended problem, in which there is an additional dependency on the correct choice of the right-hand side for constructing the field over the fracture.

A different approach to deal with crack propagation problems and avoid the smeared interface is the use of level-sets and the eXtended finite element method (XFEM) to resolve jumps in the deformation across the crack \cite{GMB02,SC03,SCBM08,FB11}. However, these approaches have in common that the crack is considered to be a lower dimensional manifold, whereas we aim to consider an open crack filled with a fluid. We also observe that there have been some approaches which couple both the phase-field and level-set/extended finite element (XFEM) approaches \cite{YNK20,LWW17,GSF17}, or alternatively, the phase-field method coupled with the numerical manifold method for the explicit construction of the lower dimensional crack. Finally, we also note that some work has gone into a pure level-set algorithm to determine the front of an open hydraulically driven fracture~\cite{PD08}. Conceptionally, level-set methods are close to phase-field methods. We prefer the latter because thermodynamic laws can be established for certain phase-field fracture models \cite{MieWelHof10a,MieWelHof10b}.

In this work, we aim to overcome the difficulties presented by the smeared interface in the context of multi-physics crack propagation. 
To this end, we couple advanced meshing software
with developments in the previously mentioned phase-field fracture methods. Our main objective is to design a high-accuracy framework in which the interface 
is computed with the help of the phase-field function. To realise the accurate interface reconstruction, we design multiple approaches to recover the interface from the phase-field approach to allow an exact description of the physics on the interface between the fluid-filled fracture and the surrounding solid. Based on this description of the interface, we then construct new meshes where the interface is aligned with mesh edges. Such alignments 
with the mesh are very well-known from interface-tracking approaches such as the arbitrary Lagrangian-Eulerian (ALE) framework \cite{HuLiZi81,DoGiuHa82} and $r$-adaptivity \cite{BuHuaRu09}. However, we also note that ALE is a more general 
concept, which is not limited to interface problems, but is 
a technique that combines Lagrangian and Eulerian coordinate systems; see, for example, 
the overview chapter \cite{DoHuePoRo04}.
A related work is \cite{Wi16_fsi_pff} in which 
the same type of phase-field fracture models are combined with 
the same type of fluid-structure interaction. However, therein, the idea was to treat 
certain interfaces with phase-field (interface-capturing) and others with ALE (interface-tracking). This is a specific situation of the general framework for interface-capturing 
and interface-tracking approaches presented earlier for finite element computations 
in \cite{Tez06}.
To the best of our knowledge, the proposed reconstruction technique from interface-capturing to interface-tracking for phase-field fracture applications is the first work in this direction.

The outline of this paper is as follows: 
In Section \ref{sec_pff}, we introduce 
the phase-field fracture framework for pressurised elastic solids.
Next, in Section \ref{sec.crack_reconstruct}, we present the details of the design of our approaches to recover the sharp interfaces based on the diffusive phase-field.
With the sharp interface at hand, we can model and couple different equations.
As an example in this work, we consider, in Section \ref{sec_fsi}, a Stokes 
flow in the fracture and coupled to the surrounding elastic medium, yielding a fluid-structure interaction problem. As we have resolved the interface between the fluid and the solid phases with our mesh and we wish to continue to track the interface, we couple the two phases via the ALE method. Finally, in Section \ref{sec_tests}, we conduct various numerical tests to demonstrate the feasibility of our algorithm. First, we focus on the interface reconstruction only while studying 
feasibility and accuracy. In the second example, we couple to a Stokes flow. 
In the third numerical test, the Stokes flow is coupled to the elastic solid, 
yielding a fluid-structure interaction problem.
Finally, in the fourth example, we extend to two orthogonal fracture, which is not challenging for the phase-field method, but new and difficult for our proposed reconstructions, since the crack opening displacement evaluation and geometry reconstruction become more involved.
Our work is summarised in Section \ref{sec_conclusions}.

\section{Modelling Fractured Solids With a Phase-Field Method}
\label{sec_pff}

The first component of our approach, is the approximation of the crack dynamics
through a phase-field approach. Based on the phase-field approximation, we
will then reconstruct the geometry of the multi-physics problem.
In this section, we therefore introduce the basic notation and  
the underlying equations for modelling fractured solids with 
a phase-field method, as well as the quantities of interest that define
the aperture of an opening crack. This in turn will then provide the
information needed to reconstruct the opening fracture for the multi-physics
simulation.

\subsection{Notation, Domains, and Regularising Interfaces}
In the following, let $\O \subset \RR^2$ the total domain under consideration.
Let  $\mathcal{C}\subset\O$ denote the fracture in our domain
and $\SO \subset \O$ is the intact domain.
For simplicity, we will assume homogeneous Dirichlet 
conditions on the outer boundary $\partial\O$.
In a phase-field approach, the fracture $\mathcal{C}$
is approximated by $\CR\subset\RR^2$ with the help of an elliptic
(Ambrosio-Tortorelli) functional \cite{AmTo90,AmTo92}.
For fracture formulations posed in a variational setting, this was first
proposed in \cite{BourFraMar00}. 
The inner fracture boundary is denoted by $\partial\CR$. 
We emphasise that the domains $\SO,\CR$, and the boundary $\partial\CR$ depend
on the choice of the so-called \emph{phase-field regularisation parameter}
$\ep>0$. Details of this parameter are presented below.
Finally, we denote the $L^2(\O)$ scalar product with $(\cdot, \cdot)$.

\subsection{Weak Formulations for the Phase-Field Crack Propagation Problem}
Weak formulations are stated in an incremental (i.e., time-discretised) formulation,
based on a quasi-static variational fracture model \cite{FraMar98,BourFraMar00}
with an extension to pressurised fracture presented in \cite{MiWheWi13a,MWW19}.
To present this, we  first introduce a classical formulation and then state 
a linearised, regularised second formulation that we use in our implementation.

The problem is formulated in two unknowns: a vector-valued displacement field
$\ub$ and a scalar-valued phase-field function $\varphi$. The phase-field function indicates the presence of a crack by taking the value $\varphi=1$ in the intact domain $\SO$, the value $\varphi=0$ inside the crack $\CR$ and a smooth transition between the two in a region of width $\ep>0$ around the interface between the open crack and the intact domain, 
denoted 
now by $\SO$ and $\CR$. Later, these domains will become $\SO$ and $\FL$, respectively, 
when a fluid is described in the crack region.
Consequently, the phase-field 
is subject to a crack irreversibility constraint $\partial_t \varphi \leq 0$.
To derive our phase-field model, the continuous irreversibility constraint 
is approximated  with a difference quotient by
\[
\varphi \leq \varphi^{old}.
\]
Here, $\varphi^{old}$ will later denote the solution at the previous time 
step $\varphi^{n-1}$ and the current solution
$\varphi:=\varphi^n:=\varphi(t_n)$ at the time point $t_n$ for $n=1,\ldots,N$. 
Now, let $\Vb\coloneqq [H^1_0(\O)]^2$ and 
\[
K\coloneqq \{w\in H^1(\O) |\, w\leq \varphi^{old} \leq 1 \text{ a.e. on }\O\},
\] 
be the solution sets, namely a function space and a convex set, respectively.
For later purposes after regularising $\varphi \leq \varphi^{old}$ by penalisation,
we also define
$W\coloneqq H^1(\O)$. The resulting system is a coupled variational
inequality system (CVIS) \cite{Wic20} and reads
\begin{form}
\label{form_1}
Let $p\in L^{\infty}(\O)$, Dirichlet boundary data $\ub_D$ on $\partial\Omega$, 
and the initial condition $\varphi(0)\coloneqq\varphi_0$ be given.
For the loading steps $n=1,2,3,\dots, N$, we compute:
Find $(\ub,\varphi)\coloneqq(\ub^n,\varphi^n) \in \{\ub_D + \Vb\} \times K$ such that
\begin{align*}\label{d}
 \Bigl(g(\varphi) \;\sigb_s(\ub), \eb( {\wb})\Bigr)
    + ({\varphi}^{2} p, \nabla\cdot  {\wb}) &= 0 \quad \forall \wb\in \Vb,\\
  \begin{multlined}[b]
    (1-\kappa) ({\varphi} \;\sigb_s(\ub):\eb(\ub), \psi {-\varphi}) 
    +  2 ({\varphi}\;  p\; \nabla\cdot  \ub,\psi{-\varphi})\\
    + G_c  \Bigl( -\frac{1}{\ep} (1-\varphi,\psi{-\varphi}) + \ep (\nabla
    \varphi, \nabla (\psi - {\varphi}))   \Bigr)  
  \end{multlined}
    &\geq  0 \quad \forall \psi \in K\cap L^{\infty}(\O).
\end{align*}
Therein, we have first the degradation function 
\[
g(\varphi)\coloneqq (1-\kappa) {\varphi}^2 + \kappa,
\]
the bulk regularisation parameter $\kappa>0$,
the phase-field regularisation parameter $\ep>0$ (linked to the spatial 
mesh size $\ep > h$ after discretisation), the Cauchy stress tensor 
\[
\sigb_s = 2 \mu \eb(\ub) + \lambda tr(\eb(\ub)) I,
\]
with the Lam\'e parameters $\mu,\lambda>0$, the identity matrix 
$I\in\mathbb{R}^{2\times 2}$
and the linearised strain tensor 
\[
\eb(\ub) = \frac{1}{2} (\nabla \ub + \nabla \ub^T).
\]

\end{form}
This system does not explicitly contain time-derivatives. Rather, 
the time $t$ might enter through time-dependent boundary conditions, e.g.,
$\ub_D=\ub_D(t)=\gb(t)$ on $\partial\O$ with a prescribed boundary function 
$g(t)$ of Dirichlet-type or  
through time-dependent right-hand side forces, e.g., a 
time-dependent pressure force $p\coloneqq p(t)$.
In this context, the latter is of interest.
Due to the quasi-static nature of this problem formulation, we derive 
the time-discretised formulation with some further approximations. 
Our first approximation relaxes the non-linear behaviour in the
first term $g(\varphi) \;\sigb_s(\ub)$ of the displacement equation by using 
\[
\varphi \approx \varphi^{n-1}, 
\]
yielding $g(\varphi^{n-1})$.
This idea is based on the 
extrapolation introduced in \cite{HWW15} and is numerically justified 
specifically for slowly growing fractures \cite[Chapter 6]{Wic20}, 
while counter examples for fast-growing fractures were found in \cite{Wic17}.
Since we are mainly interested in crack width variations in this paper, and less 
in variations in the length, our approximation is numerically justified. 
The second approximation is related to the inequality constraint. In this 
work, we relax the constraint by simple penalisation \cite{MWT15} (see also 
\cite[Chapter 5]{Wic20}), i.e., 
\[
\varphi \leq \varphi^{old} \quad\rightarrow\quad \gamma(\varphi - \varphi^{n-1})^+.
\]
Here, $(x)^+ = x$ for $x>0$ and $(x)^+ = 0$ for $x\leq 0$, 
and where $\gamma>0$ is a penalisation parameter.
We then arrive at the regularised scheme
\begin{form}\label{form_2}
Let $p\in L^{\infty}(\O)$ and the initial condition $\varphi(0)\coloneqq\varphi_0$ be given.
For the loading steps $n=1,2,3,\dots, N$, we compute:
Find $(\ub,\varphi)\coloneqq(\ub^n,\varphi^n) \in \Vb \times W$ such that
\begin{align*}
  \Bigl(g(\varphi^{n-1})\;\sigma(\ub), \eb({\wb})\Bigr)
    +({\varphi^{n-1}}^{2} p, \nabla\cdot{\wb}) &= 0
  \quad \forall \wb\in \Vb,\\
  \begin{multlined}[b]
    (1-\kappa) ({\varphi} \;\sigma(\ub):\eb(\ub), \psi) 
      +  2 ({\varphi}\;  p\; \nabla\cdot  \ub,\psi)\\
    +  G_c  \Bigl( -\frac{1}{\ep} (1-\varphi,\psi) 
      + \ep (\nabla\varphi, \nabla\psi) \Bigr) 
      + (\gamma(\varphi - \varphi^{n-1})^+,\psi)
  \end{multlined}
   &=0\quad\forall \psi\in W.
\end{align*}
\end{form}
\Cref{form_2} will then be the system of equations we use to compute the phase-field fracture.
We notice that $N=5$ is used in this work, because stationary fractures are our main interest, with the main goal to construct from the phase-field interface-capturing 
technique, a subsequent interface-tracking representation. The extension to propagating 
fractures with $N\gg 1$ is left for future work.

\subsection{Crack aperture and volume}
We briefly discuss central quantities of interests which are easily recoded from the phase-field model and provide us with quantitative details of the crack geometry. These quantities will then also be central to the sharp interface reconstruction of the opening crack interface.

The \emph{crack opening displacement} (COD), or aperture of the crack, can be computed from the phase-field by
\begin{equation*}
 \COD(\xb) = \llbracket \ub\cdot\nb \rrbracket
 \simeq \int_{\ell^{\xb,\bm{v}}} \ub(\xb) \cdot \nabla \varphi(\xb) \dif s,
\end{equation*}
where $\ell^{\xb,\bm{v}}$ is a line through $\xb$ along the vector $\bm{v}$ \cite{CHUKWUDOZIE2019957}, where $\bm{v}$ is in $\Omega$. The left-hand side allows under knowledge of the unit normal vector $\nb$ to compute the COD at each point $\xb\in \O$. 
We note that the normal vector $\nb$ points into the direction of $\nabla \varphi(\xb)$, 
because $\varphi$ represents the level-sets of the fracture iso-surfaces. A theoretical 
justification that the crack opening displacements can be formulated without normalization 
$\|\nabla \varphi(\xb)\|$ is provided in \cite{CHUKWUDOZIE2019957}[Section 3.2].
If we have $\O = (a,b)\times(c,d)$ and the crack $\CR$ is aligned with the $x$-axis in a Cartesian coordinate system, then the COD simplifies to
\begin{equation}\label{eqn.cod}
   \COD(\xb) = \int_c^d\ub(\xb_0, s)\cdot\nabla\varphi(\xb_0, s)\dif s,
\end{equation}
see, for example, \cite[Proposition~83]{Wic20}.

A second quantity of interest will be the \emph{total crack volume} (TCV) of the open crack.
By integrating over all lines, i.e., the entire domain, this can be computed by
\begin{equation}\label{eqn.tcv}
  \TCV(\CR) = \int_{\O}\ub\cdot\nabla\varphi\dif\xb.
\end{equation}
See~\cite[Definition~72]{Wic20}.

\section{Eulerian Crack Reconstruction}
\label{sec.crack_reconstruct}

Now that we have an approximation of the crack at hand, 
and a decomposition into $\SO\coloneqq\O\setminus\CR$ and $\FL\coloneqq\CR$,
we discuss a number of possible approaches to reconstruct the sharp interface of the opening crack.
The aim for this reconstruction is to obtain an Eulerian description of the domain geometry, which can then be used to solve the fluid-structure-interaction problem between the cracking solid and, for example, a fluid filling the crack interior.

\subsection{Level-Set Approaches}
We first consider an approach that utilises the mesh used for the phase-field simulation. This is an attractive choice, if the resulting partitioned geometry is difficult to mesh. Following this approach, a fixed grid fluid-structure interaction solver then needs to be utilised to realise the multi-physics simulation. Such approaches have the advantages over moving mesh approaches, that topology changes are feasible. However, we note that these approaches are also not as developed and thoroughly investigated as moving mesh discretisations.

To this end, we aim to construct a level-set function describing the cracked geometry. A level-set is a function $\phivar\colon\O\rightarrow\RR$, who's zero iso-surface describes the boundary of the crack, i.e.,
\begin{equation*}
  \partial\CR = \{\xb\in\O\;\vert\; \phivar(\xb) = 0\}
  \quad\text{and}\quad
  \CR = \{ \xb\in\O\;\vert\; \phivar(\xb) < 0\}.
\end{equation*}
For an opening fracture, the phase-field can be seen as a level-set function for the crack domain~\cite{LWW17}. To use this as a level-set, we need to identify the correct iso-surface, i.e., determine $c_{ls}$ such that
\begin{equation*}
  \phivar = \varphi - c_{ls}. 
\end{equation*}
In \cite{LWW17}, $c_{ls}$ was chosen as $0.1$, while in \cite{YNK20}, the choice was refined to $c_{ls} = \frac{\sqrt{5} - 1}{2}\approx 0.618$ by analysing a one dimensional problem. However, this does not provide an Eulerian description of the crack interface and information from the displacement is needed. The crack opening displacement is given by the normal displacement on the level-set. On a line $\ell$ perpendicular to the centreline, described by a level-set $\phivar_\ell$, the COD is
\begin{equation}\label{eqn.cod.lset}
  \COD(\ell) = \meas_{0}\left(\{\ub\cdot\nb_{\phivar}\;\vert\; \phivar(\xb) = 0 \text{ and } \phivar_{\ell}(\xb)=0\}\right),
\end{equation}
where $\nb_\phivar = \nabla\phivar / \nrm{\nabla\phivar}$ 
is the normal vector on the level-set boundary. As a result, we see that the iso-surfaces of the phase-field function are not feasible as an Eulerian level-set description of the crack boundary.

\subsubsection{Explicit Level-Set Construction}
\label{sec.crack_reconstruct:subsec.lset:explicit}
Based on the knowledge of the centreline and the crack opening displacements, we obtain a set of points on the Eulerian cracks interface $\{(\xb_1^i, \xb_2^i) \}$. Based on these points, we can then create a set of lines $\phivar_i$ that connect two adjacent points $(\xb_1^{i-1}, \xb_2^{i-1})$, $(\xb_1^i, \xb_2^i)$ with straight lines for $i=0, \dots, n$ with $(\xb_1^{-1}, \xb_2^{-1})=(\xb_1^{n}, \xb_2^{n})$. Extending each of these line segments to lines in the whole domain, the crack interior becomes a geometry described by multiple level-set functions, from which a single level-set can then be reconstructed~\cite{BCH14}. For example, if we are interested in the domain, where two level-sets are negative, we can use the single level-set $\phivar = \min\{\phivar_1, \phivar_2\}$. The resulting function can then be cast into a single piece-wise linear function on the mesh by an appropriate piece-wise linear interpolation. A sketch of this construction idea can be seen on the left of \Cref{fig.geometry-approaches}.

\subsubsection{Level-Set Transport}
\label{sec.crack_reconstruct:subsec.lset:transport}
A different approach to construct an Eulerian level-set description of the domain is to transport the phase-field level-set $ \varphi - c_{ls}$ to the Eulerian level-set. Assuming that the crack 
is aligned with the $\xb_1$-axis, then we know that each point $\xb=(\xb_1, \xb_2)$ on zero-line of the phase-field level-set must be mapped to $(\xb_1, \ub(\xb)\cdot\nb_{\phivar} / 2)$. To move the level-set points $\{\xb\}$ to $\{(\xb_1, \ub(\xb)\cdot\nb_{\phivar} / 2)\}$, we transport the level-set along a velocity field $\bm{\beta}$. This is realised by solving the transport problem
\begin{equation*}
  \partial_t\phivar + \bm{\beta}\cdot\nabla\phivar = 0.
\end{equation*}
This can then be numerically treated, for example by the Streamline-Upwind-Petrov-Galerkin (SUPG) method~\Cite{Bur10a} or flux-corrected transport (FCT) \cite{Kuzmin20092517}. For this approach to work in the setting, we 
need to construct an appropriate transport field $\beta$. Since $\xb \mapsto (\xb_1, \ub(\xb)\cdot\nb_{\phivar} / 2)$, we construct
\begin{equation*}
  \beta = \begin{pmatrix}
    0\\
    \ub(\xb)\cdot\nb_{\phivar} / 2 - \xb_2,
  \end{pmatrix}
\end{equation*}
on the zero iso-surface of $\phivar$ and extend this harmonically into the volume with zero Dirichlet boundary conditions on the outer boundary. We then consider the artificial time-interval $[0, t]=[0,1]$ over which to solve the transport problem. As the transport field is not constant in space, we also transport the components of the velocity field, so that the transport along the lines $\nb_\phivar$ is constant in each time-step and that we reach the correct line after time $t=1$. A sketch of this construction idea can also be seen in the centre of \Cref{fig.geometry-approaches}.

\begin{remark}
  All the level-set approaches are only feasible if the aperture of the crack is larger than the mesh size $h$. This is due to the fact that most unfitted finite element codes construct the unfitted quadrature rules under the assumption, that each element is cut by the level-set function at most once. Furthermore, the finite element space with which we solve the fluid problem inside the crack needs to be sufficiently large to have the necessary approximation properties.
\end{remark}

\begin{remark}
Each point $\xb$ on zero-line of the phase-field level-set must be mapped to $(\xb_1, \ub(\xb)\cdot\nb_{\phivar} / 2)$. This could in theory also be realised by a deformation of the original mesh, so that the level-set $\phivar = \varphi - c_{ls}$ is an Eulerian description of the crack on the deformed mesh in world coordinates. Unfortunately, the deformation that would move the zero iso-line accordingly is too large near the tip edges, such that the resulting mesh becomes degenerated.
\end{remark}

\subsection{Explicit Mesh Construction}
\label{sec.crack_reconstruct:subsec.expl_geom}
Based on the explicit knowledge of the points on the crack interface as in \Cref{sec.crack_reconstruct:subsec.lset:explicit}, we can also generate a geometry description which we can then re-mesh to create a mesh fitted to the crack/fluid and solid domains.

\begin{figure}
  \centering
  \includegraphics[width=5.2cm]{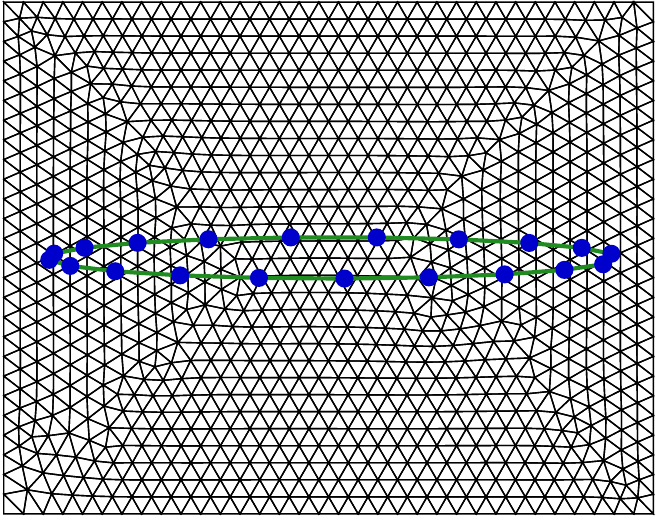}
  \includegraphics[width=5.2cm]{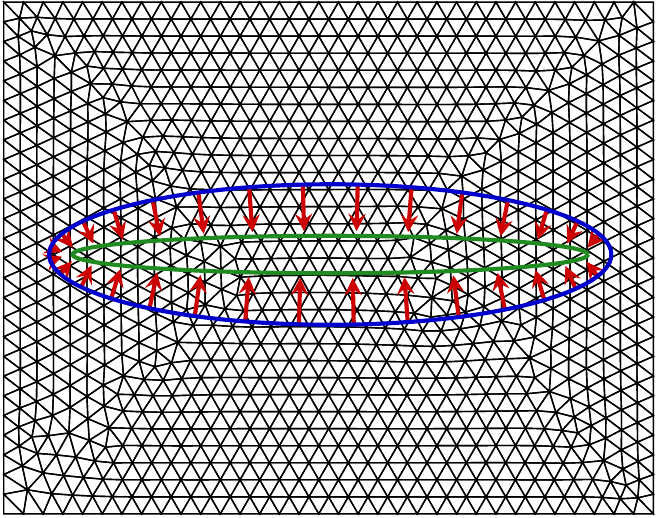}
  \includegraphics[width=5.2cm]{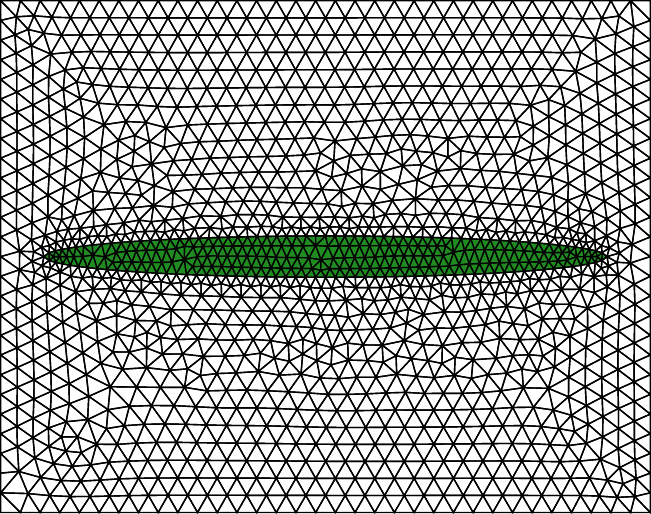}
  \caption{Level set and mesh construction approaches. Left: Explicit Level set construction based on a set of points; Center: Level-Set Transport from the blue to the green iso-surface; Right: Explicit mesh construction with elements inside the crack marked green.}
  \label{fig.geometry-approaches}
\end{figure}

To realise this, we can construct a piecewise linear description of the crack boundary, similar to the explicit level-set construction in \Cref{sec.crack_reconstruct:subsec.lset:explicit}. Alternatively, we can also realise a curved crack boundary. This can be achieved in \texttt{netgen}~\cite{Sch97} by using a description with rational splines of second order, such that the tangent of the resulting curve is continuous at all supporting points on the boundary. A sketch of the mesh resulting from this approach can be seen in the right of \Cref{fig.geometry-approaches}.

Depending on the number of supporting points and shape of the crack, the generation of the mesh fitted to the crack geometry can be a challenging task. However, once we have constructed such a mesh, then the well studied arbitrary Lagrangian-Eulerian approach (see \Cref{sec_fsi} below) becomes feasible for solving a coupled fluid-structure interaction problem between a fluid in the crack and the elastic material cracking; for example, by taking the Eulerian domain as the reference domain $\Oref$.

\section{Coupling Fracture Stokes Flow to the Surrounding Elastic Medium}
\label{sec_fsi}

With the mesh containing the resolved interface between the open crack and the
surrounding solid at hand, we can now describe the fluid-structure interaction
problem between the fluid-filled crack and the solid domain surrounding it.
With our explicit interface reconstructions of the fracture 
surface from \Cref{sec.crack_reconstruct:subsec.expl_geom}, the flow problem is coupled via interface-tracking to the surrounding
solid, and we arrive at a classical fluid-structure interaction model.
In order to couple flow and solids, we discuss the
arbitrary Lagrangian-Eulerian approach below. Finally, at the end of this
section, we gather all previous algorithms and design a final overall scheme of our 
high-accuracy phase-field interface-reconstruction framework.

\subsection{Stationary Fluid-Structure Interaction}
In this section, we model fluid-structure interaction in 
arbitrary Lagrangian-Eulerian coordinates using variational 
monolithic coupling in a reference configuration \cite{HrTu06a,Du07,Wi11_phd,Ri17_fsi}.
For simplicity, we shall consider the linear Stokes equations as the fluid model here.

Consider a domain $\O\subset\RR^d$ divided into a $d$-dimensional fluid domain $\FL$, a $d$-dimensional solid domain $\SO$ and a $d-1$-dimensional interface $\IN$ between the two, such that $\O=\FL\dot\cup\IN\dot\cup\SO$. Furthermore, let $\hat{\O}, \FLref, \SOref$ and $\INref$ be the corresponding domains in a reference configuration. In our setting the fluid domain is the interior of the crack $\FL=\CR$, the solid is the untracked medium $\SO$ and the interface is the crack boundary $\IN=\partial\CR$.

Using the reference domains $\SOref$ and $\FLref$ leads to
the well-established ALE coordinates \cite{HuLiZi81,DoGiuHa82}. To obtain a monolithic
formulation we need to specify the transformation $\hcalA_f$ in the
fluid-domain. On the interface $\INref$ 
this transformation is given by the structure displacement:
\[
\hcalA_f(\xbref,t)\big|_{\INref} = 
\xbref+\ubref_s(\xbref,t)\big|_{\INref}. 
\]
On the outer boundary of the fluid domain
$\partial\FLref\setminus\INref$ it holds $\hcalA_f=\text{id}$.
Inside
$\FLref$ the transformation should be as smooth and regular as
possible, but apart from that it is arbitrary. Thus we harmonically
extend $\ubref_s|_{\SOref}$ to the fluid domain $\FLref$ and
define $\hcalA_f\coloneqq \text{id}+\ubref$ on $\FLref$,
where $id(\xbref) = \xbref$ in $\hcalA_f\coloneqq \text{id}+\ubref
:= \hcalA_f(\xbref,t)\coloneqq \text{id}(\xbref,t)+\ubref(\xbref,t)$
such that
\[
(\hat \nabla \ubref_f,\hat \nabla\psiref)_{\FLref} 
 = 0,\quad
\ubref_f=\ubref_s\text{ on }\INref,\quad
\ubref_f=0\text{ on }\partial\FLref\setminus\INref.
\]
Consequently, we define a continuous variable $\ubref$ on  all $\Omega$
defining the deformation in $\SOref$ and supporting the
transformation in $\FLref$. By skipping the subscripts and
since the definition of $\hcalA_f$ coincides with the definition of
the solid transformation $\hcalA_s$, we define on
all $\Oref$:
\begin{equation*}
  \hcalA(\xbref, t)\coloneqq \xbref + \ubref(\xbref, t),\quad
  \hat F(\xbref, t) \coloneqq \hat\nabla \hcalA= I+\hat\nabla \ubref(\xbref, t),\quad
  \hat J\coloneqq \text{det}(\hat F).
\end{equation*}

With this at hand, the weak formulation of the stationary fluid-structure interaction problem is given by \cite{RiWi10}:
\begin{form}[Stationary fluid-structure interaction]
  \label{fsi:ale:stationary}
  Let $\Vbref$ be a subspace of $\bm{H}^1(\Oref)$ with trace zero on
  $\hat\Gamma^D\coloneqq \hat\Gamma_f^D\cup\hat\Gamma_s^D$ and $\hat L\coloneqq L^2(\Oref)/\mathbb{R}$. Find $\vbref\in \Vbref$, $\ubref\in
  \Vbref$ and $\hat p\in \hat L$, such
  \begin{subequations}\label{eqn.ale-fsi}
    \begin{align}
(\hat J\sigbref_f\hat F^{-T},\hat \nabla\phiref)_{\FLref}
      + (\hat J\sigbref_s\hat F^{-T},\hat \nabla\phiref)_{\SOref}
&= (\rho_f \hat J\fbref,\phiref)_{\FLref} &&\forall\phiref\in\Vbref,\\
- (\vbref,\psiref)_{\SOref} + 
      (\alpha_u \hat \nabla \ubref,\hat \nabla\psiref)_{\FLref}
&=0&&\forall\psiref\in\Vbref,\\
(\widehat{\diver}\,(\hat J\hat F^{-1}\vbref_f),\xiref)_{\FLref} 
&=0&&\forall\xiref\in \hat L,
    \end{align}
  \end{subequations}  
  with a right-hand side fluid force $\fbref\in L^2(\FLref)$ and the 
harmonic mesh extension parameter $\alpha_u>0$. 
Finally, the Cauchy stress tensor in the solid is
defined in \Cref{form_1} and we use
$\hat J\sigbref_s\hat F^{-T} \coloneqq  \sigb_s$.
The ALE fluid Cauchy stress tensor $\sigma_f$ is given by 
\[
  \sigbref_f \coloneqq  -\hat p_fI +\rho_f\nu_f(\hat\nabla \vbref_f \hat F^{-1}
  + \hat F^{-T}\hat\nabla \vbref_f^T),
  \]
  with the kinematic viscosity $\nu_f$ and the fluid's density $\rho_f$.
\end{form}

\subsection{Final Algorithm}
With the derivations of the previous sections, we formulate the following scheme:
\begin{Algorithm}\label{alg.full}
~\\\vspace{-10pt}
\begin{enumerate}
  \item Given some pressure $p$, compute the phase-field approximation of the crack $(\ub,\varphi)$ using \Cref{form_2}.
  \item Reconstruct the sharp crack interface $\partial\CR$ using the approach presented in \Cref{sec.crack_reconstruct:subsec.expl_geom} and re-mesh the resulting geometry.
  \item Assign reference domains: $\Oref=\O, \FLref=\CR$ and $\SOref=\O\setminus\CR$, where $\O$ is the domain used to compute the phase-field using \Cref{form_2} in Step 1.
  \item Given some right-hand side fluid-force $\fbref$ in $\FLref$, compute the fluid-structure interaction problem using the ALE approach in \Cref{fsi:ale:stationary}.
\end{enumerate}
\end{Algorithm}

\begin{remark}
After the FSI step No.\@{} 4, the next step would be to take the Stokes pressure $p$ 
and to go to Step 1, which would result into an iterative loop in which 
phase-field and interfaces are successively corrected. 
The practical realisation and computational analysis of this loop is left for future work. 
This opens the way to implement a time-stepping for time-dependent
situations with $t_n, n=1,\ldots, N$ 
with propagating fractures and non-stationary fluid-structure interaction.
\end{remark}

\section{Numerical Tests}
\label{sec_tests}
In this section, we conduct several numerical experiments. First, we consider 
Sneddon's test \cite{SneddLow69}, which is nowadays considered as a benchmark problem \cite{Schroeetal20}.
Therein, the pressure is a given quantity.
In the second numerical example, we further investigate the quality of our domain reconstruction approaches 
by solving the Stokes equations inside the fracture.
In the third and fourth numerical tests, we consider the full \Cref{alg.full}, i.e., a stationary fluid-structure interaction problem based on the domain reconstructed from the phase-field approximation.

Our examples are implemented using \texttt{Netgen/NGSolve} \cite{Sch97,Sch14} together with the add-on \texttt{ngsxfem} \cite{LHPvW21} for unfitted finite elements.

\subsection{Sneddon's Test}
\label{sec_test_sneddon}
Our first example is based on Sneddon's theoretical calculations in \cite{Sne46,SneddLow69}. Specifically, we consider the two-dimensional case with constant pressure $p$ acting on the fracture boundary.
In Sneddon's test, usually, the domain and all parameters 
are provided in dimensionless values, which we follow in this work as well.

\subsubsection{Configuration}
The domain is $\Omega=(0,4)^2$. The problem is stationary, as the data driving the crack is constant. For the boundary conditions, we have
\begin{align}
  \ub &= 0\quad\text{on } \partial\O,\label{eqn.senddon_bc1}\\
  \epsilon\partial_{\nb}\varphi &= 0\quad\text{on }\partial\O\label{eqn.senddon_bc2}.
\end{align}
For the initial condition, we set the phase-field as
\begin{equation}\label{eqn.sneddon_initial_pf}
  \varphi(\xb) = 
  \begin{cases}
    0 &\xb\in\CR_0 = (1.8, 2.2)\times(2-h, 2+h)\\
    1 &\xb\in\O\setminus\CR_0,
  \end{cases}
\end{equation}
with the local mesh size $h$, i.e., a crack of length $2l_0 = 0.4$ in the centre of the domain, parallel to the $x$-axis.

The mechanical parameters are  Young's modulus and Poisson's ratio which we set to be $E_s=10^{5}$ and $\nu_s=0.35$. The applied pressure is $p=4.5 \times 10^{3}$ and the critical energy release rate is chosen as $G_c=500$.

\subsubsection{Discretisation}
We use an unstructured triangular mesh of the domain that resolves $\CR_0 = (2-l_0, 2-h)\times(2+l_0,2+h)$. The mesh is constructed such that the simplicial at the crack are smaller than those at the outer boundary the domain by a factor 100. On this mesh, we use piecewise linear, continuous, finite elements for both the displacement and phase-field. The penalisation parameter is chosen as $\gamma = 100 h^{-2}$ and the phase-field regularisations parameter is set to $\kappa=10^{-10}$ and $\epsilon=0.5\sqrt{h}$. We iterate the phase-field problem for a total of five pseudo time-steps to arrive at the stationary solution.

We compute Sneddon's test over a series of meshes with $h=h_0 \cdot 2^{-l}$, $l=0, \dots,5$ and $h_0=0.02$. The resulting crack is then reconstructed from the resulting phase-field solution using the approaches discussed in \Cref{sec.crack_reconstruct}. We will refer to the level-set construction via the crack opening displacement and straight line segments as the \emph{explicit level-set} construction, c.f., \Cref{sec.crack_reconstruct:subsec.lset:explicit}. We call the approach to construct the level-set recovered from the phase-field and level-set transport as the \emph{transport level-set} approach, c.f., \Cref{sec.crack_reconstruct:subsec.lset:transport}. Finally, we refer to the approach of constructing a new  mesh the crack geometry from the spline approximation based on the crack opening displacements as the \emph{explicit mesh} approach, c.f.~\Cref{sec.crack_reconstruct:subsec.expl_geom}.

\subsubsection{Quantities of Interest}
The total crack volume of the crack resulting from Sneddon's test has the analytical expression
\begin{equation*}
  \TCV = 2\pi \frac{(1 - \nu_s^2)l_0^2p}{E_s}.
\end{equation*}
We note that this expression is only valid in an unbounded domain; see 
the derivations in \cite{SneddLow69} and a computational confirmation on the dependence 
on the domain size was carried out in \cite{HeiWi18_pamm}[Fig. 1 right].
Nevertheless, as the domain is large compared to the crack size,
we will compare this against the values computed by integrating the phase-field using formula \eqref{eqn.tcv}, and the volume of the domain constructed by the level-set or spline approximations.
For the full crack opening displacement, we have the analytical expression
\begin{equation*}
  \COD(x_0) = 4 \frac{(1 - \nu_s^2)l_0 p }{E_s}\bigg(1 - \frac{x_0^2}{l_0^2} \bigg)^{1/2}.
\end{equation*}
We will then compare this with the crack opening displacements as computed from formula \eqref{eqn.cod}, which we shall call the \emph{integration} method, and from \eqref{eqn.cod.lset}, which we will refer to as the \emph{point evaluation}.

\subsubsection{Results}
\begin{figure}
  \centering
  \includegraphics{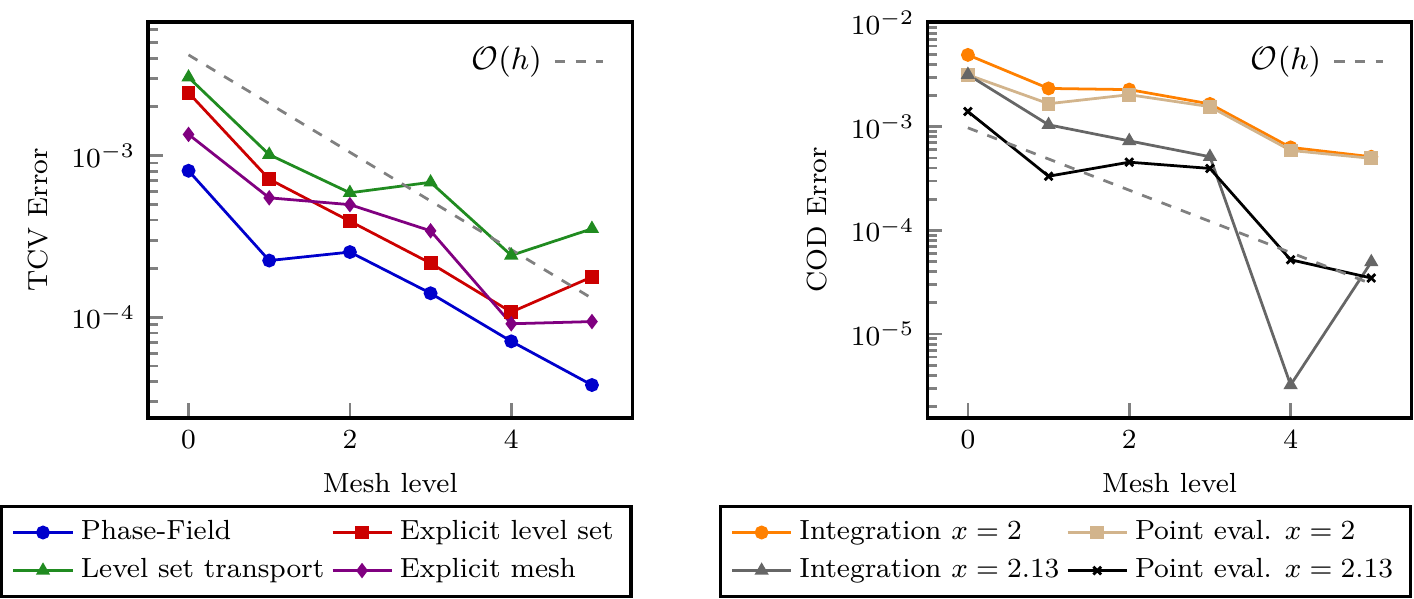}
  \caption{Example \ref{sec_test_sneddon}: Mesh convergence for Sneddon's test. Left: Total crack volume convergence computed from the phase-field integration, rational spline domain area and level-set domain areas respectively. Right: Crack opening displacement computed by integration and point evaluation of the phase-field.}
  \label{fig_sned1_tcv_cod_conv}
\end{figure}

\begin{figure}
  \centering
  \includegraphics{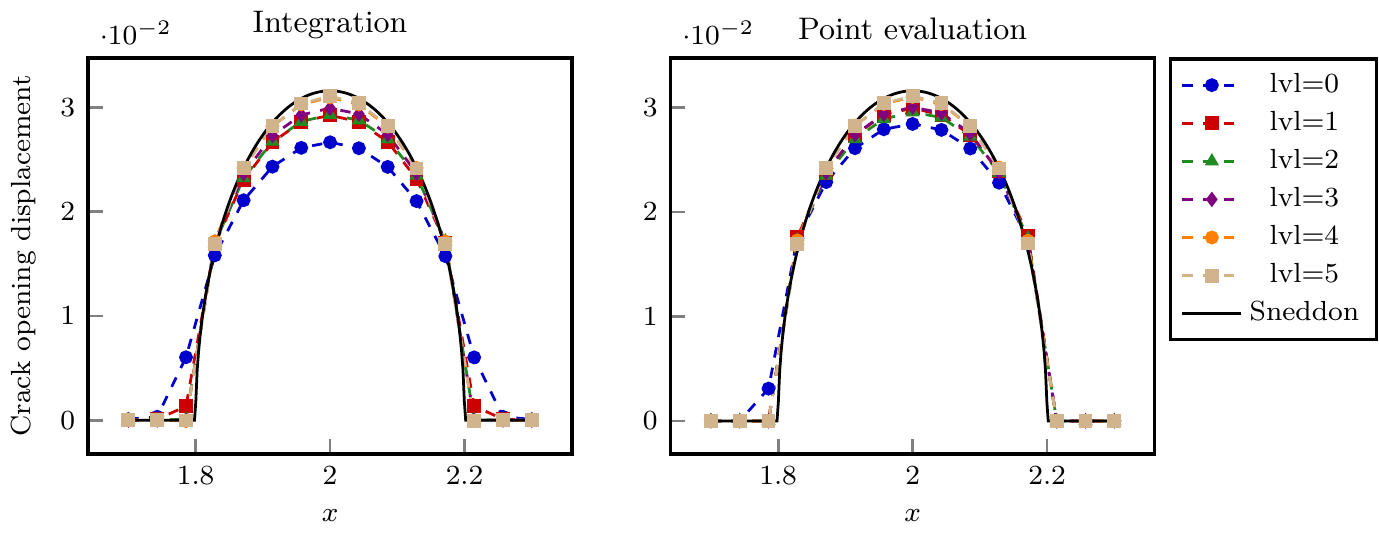}
  \caption{Example \ref{sec_test_sneddon}: Crack opening displacements at a number of points along the crack over the series of meshes considered.}
  \label{fig_sned1_cod_comp}
\end{figure}

The results of the total crack volume and crack opening displacements at the centre of the crack ($x=2$) and closer towards the tip of the crack ($x=2.13$) can be seen in \Cref{fig_sned1_tcv_cod_conv}. Furthermore, the resulting crack opening displacements as computed from the phase-field integration and point evaluation can be seen in \Cref{fig_sned1_cod_comp}.

On the left of \Cref{fig_sned1_tcv_cod_conv}, we see that the total crack volume converges towards the expected value for all four methods. Overall it appears that all methods converge linearly. The volume computed from the phase-field is the most accurate while the volume of the level-set obtained from level-set transport is the least accurate.
On the right of \Cref{fig_sned1_tcv_cod_conv}, we see that the crack opening displacement converges similarly for both approaches. In the case $x=2.13$, we see 
a numerical artefact in which the error drops too low, which is often the case 
in computational error analyses with goal functionals, as no monotone convergence 
can be expected. This behaviour is also seen in some curves in the left sub-figure of \Cref{fig_sned1_tcv_cod_conv}.
Looking at \Cref{fig_sned1_cod_comp}, we see that overall, both methods of computing the crack opening displacements capture the shape of the crack. The point evaluation appears to be more accurate on the coarsest mesh, but for finer meshes, the difference is negligible.
The main difference is that the point evaluation of the phase-field normal captures the tip of the crack more accurately on coarse meshes.

\subsection{Sneddon's Test Coupled to a Stokes Problem in the Crack}
\label{sec_test_Stokes}
For a second example, we study further the quality of the crack interface approximation and the approximation qualities of the resulting mesh on the interior of the crack. To this end, we consider the same fracture problem as above, as we know that the resulting crack is an ellipse of width $2a =0.4$ and height $2b=4 \frac{(1 - \nu_s^2)l_0 p }{E_s}$.

Following the \cite[section~3.2.6]{vWah21}, we construct an divergence free velocity field and pressure, such that the velocity conforms to the homogeneous Dirichlet boundary condition on the ellipse and the pressure has mean zero. Let $a,b$ be the semi-major and semi-minor axes of an ellipse and $\hat{\xb}$ be the centre of the ellipse. Consider the stream function $\Phi(\xb)=\sin(\nicefrac{\pi}{2} ((\xb_0 - \hat{\xb}_0)/a)^2 + ((\xb_1 - \hat{\xb}_1) / b)^2 )$. Then the vector field $\ub_\text{ex} = (\partial_{\xb_1}\Phi, -\partial_{\xb_0}\Phi)^T$ is divergence free by construction and is equal to zero on the ellipse under consideration. We further let $p_{ex}=\Phi - \nicefrac{2}{\pi}$ which has mean zero on the ellipse. The right-hand side of the Stokes problem is then set to $\fb=-\nu\Delta\ub_{ex} + \nabla p_{ex}$. We choose $\nu=10^{-4}$.

\subsubsection{Discretisation of the Stokes Problem}
To compute the Stokes problem on the level-set domain constructed from the phase-field fracture, we use unfitted finite elements known as CutFEM~\cite{BCH14}. Details of this approach for the Stokes problem is given in \cite{BH14,MLLR14}. Here we shall use the Taylor-Hood finite element pair $\PP^2/\PP^1$. Elements cut by the interface are stabilised using ghost-penalty stabilisation~\cite{Bur10} and in particular with the direct version of the ghost-penalty operator introduced in \cite{Pre18}.
In case of the re-meshed geometry, we have a fitted mesh and we use inf-sup stable Taylor-Hood elements $\PP^2/\PP^1$ on this mesh.

\subsubsection{Results}

\begin{figure}
  \centering
  \includegraphics{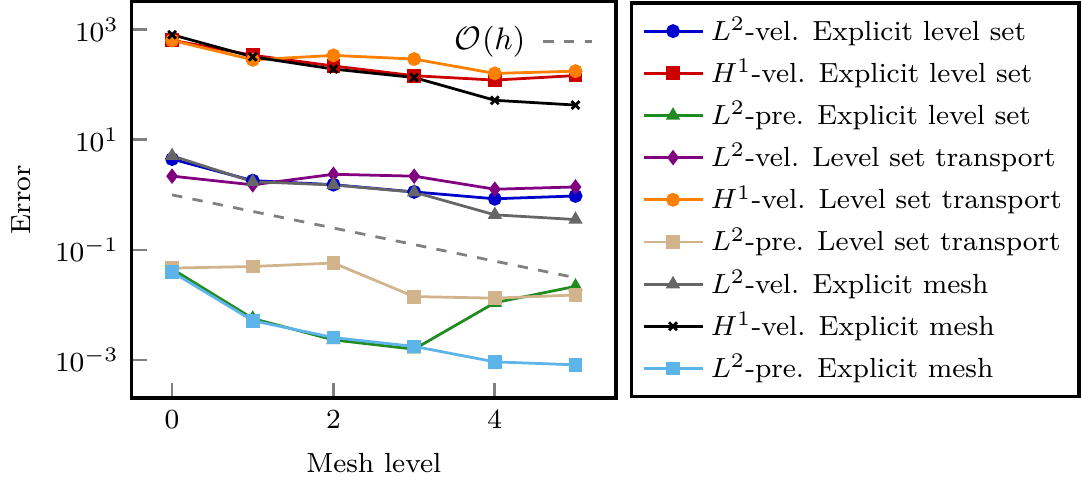}
  \caption{Example \ref{sec_test_Stokes}: Error convergence for the velocity and pressure for the Stokes problem computed on the domains resulting from Sneddon's test using the fitted and unfitted Eulerian domain approximations approaches.}
  \label{fig_stokes_on_reconstructed_domain}
\end{figure}

\begin{figure}
  \centering
  \null\hfill
  \includegraphics[width=4.5cm]{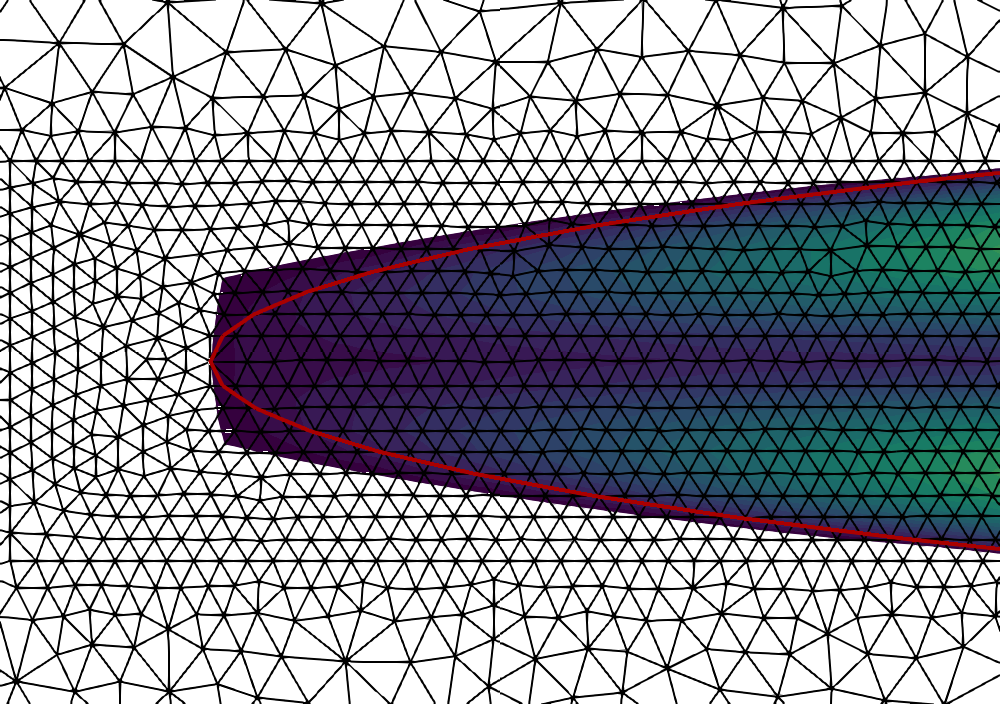}
  \hfill
  \includegraphics[width=4.5cm]{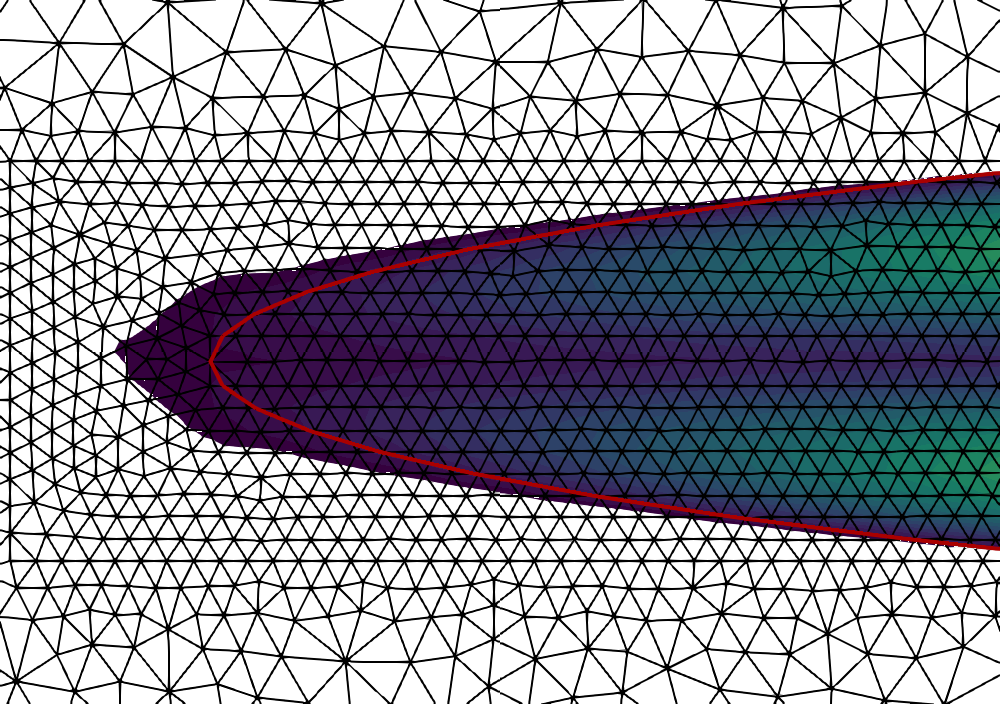}
  \hfill
  \includegraphics[width=4.5cm]{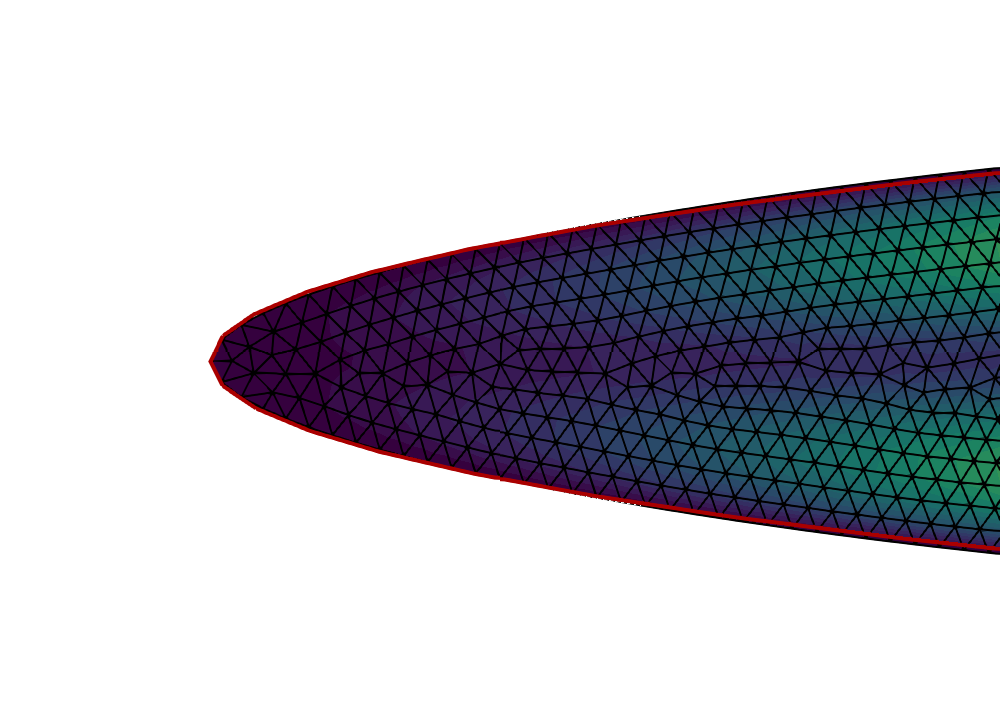}
  \hfill\null

  \vspace{8pt}
  
  \null\hfill
  \includegraphics[width=4.5cm]{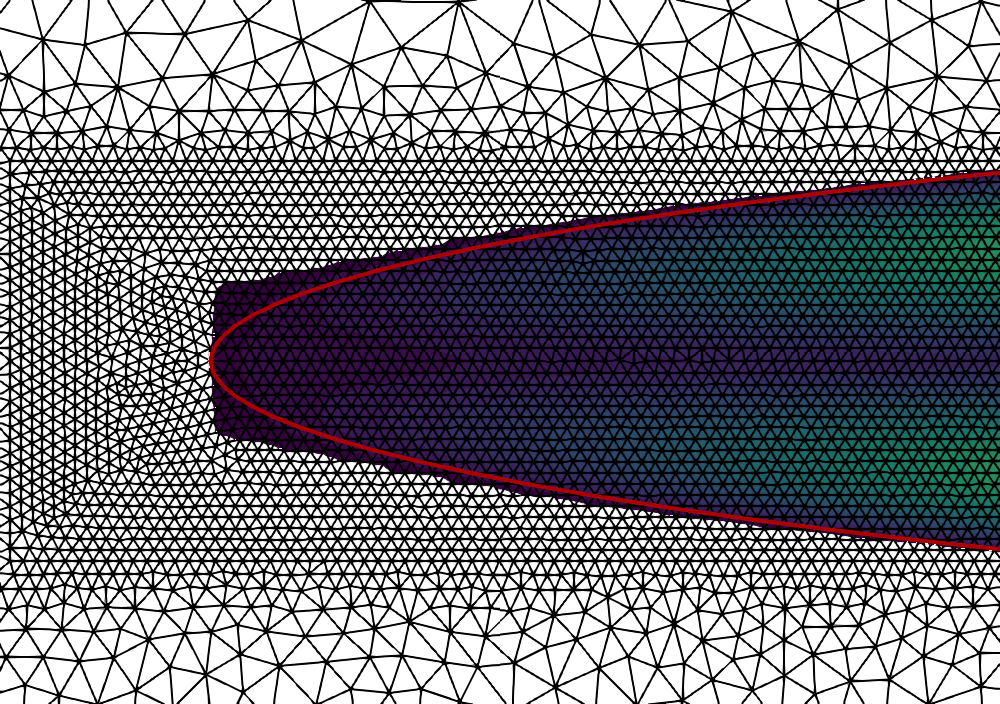}
  \hfill
  \includegraphics[width=4.5cm]{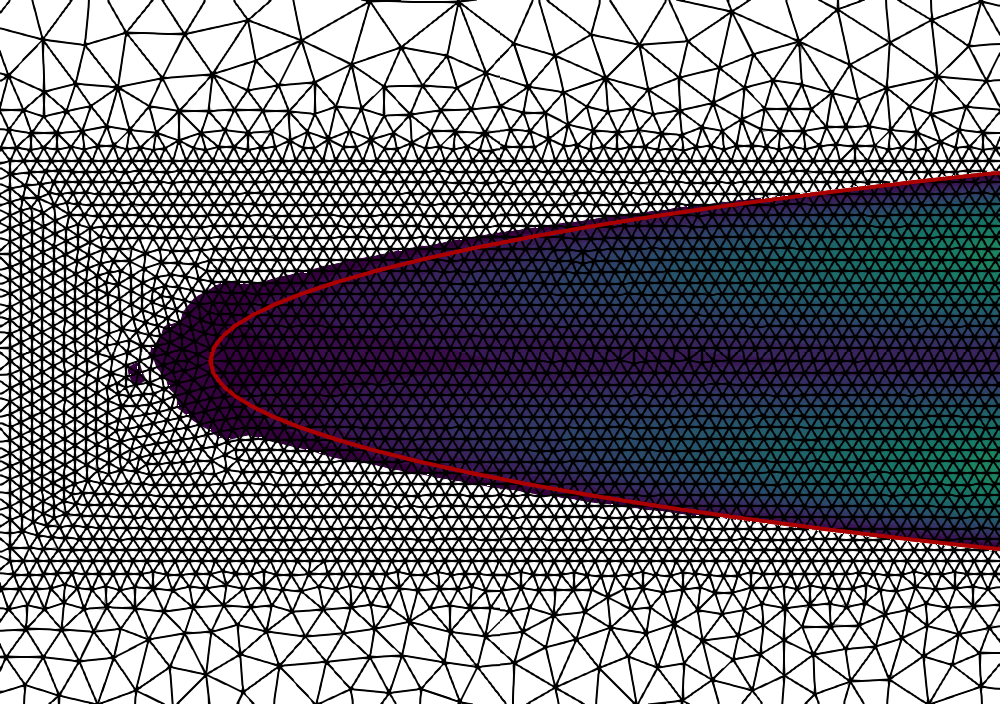}
  \hfill
  \includegraphics[width=4.5cm]{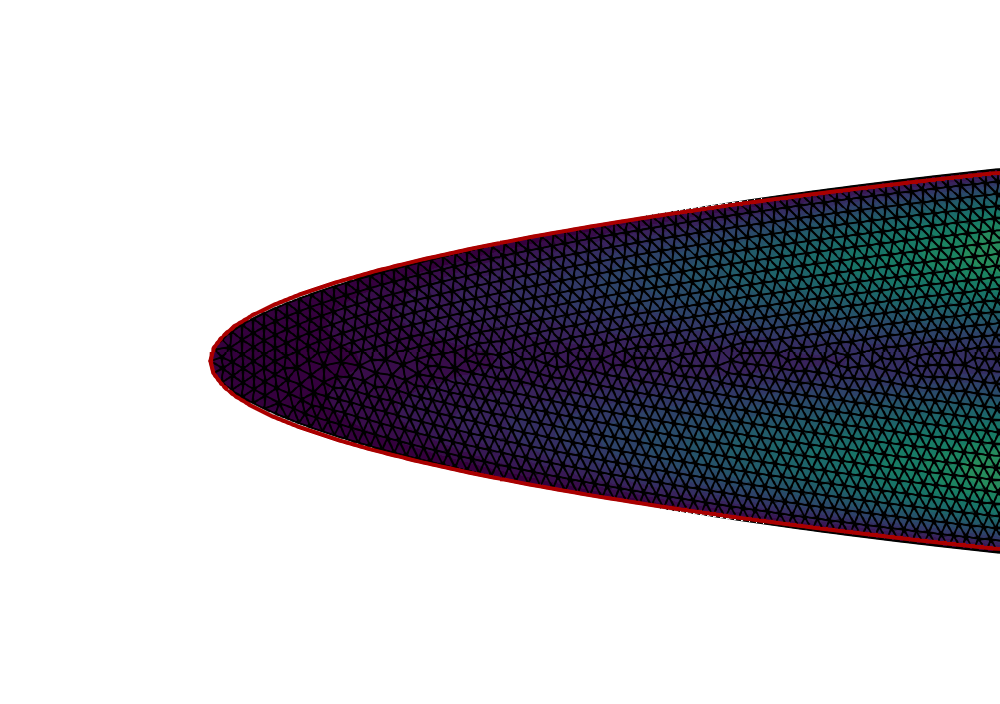}
  \hfill\null
  \caption{Example \ref{sec_test_Stokes}: Stokes velocity solution magnitude and exact ellipse (red) near the tip of the crack. Top: Mesh level 4, bottom: Mesh level 5. From left to right: Explicit level-set construction, level-set transport, rational spline construction with re-meshing.}
  \label{fig_stoke_domain_tip}
\end{figure}

The convergence results for the velocity in the $L^2$- and $H^1$-norms and of the pressure in the $L^2$-norm can be seen in \Cref{fig_stokes_on_reconstructed_domain}. Furthermore, the domains near the tip of the crack resulting from the different reconstruction approaches can be seen in \Cref{fig_stoke_domain_tip}.

In \Cref{fig_stokes_on_reconstructed_domain}, we see that we have some limited first order convergence for the resulting velocity and pressure, suggesting that the geometry error in the approximation of the crack is the limiting factor. This is consistent with our previous results in the crack aperture and volume.
We also note that the error begins to remain constant on finer meshes for the level-set approaches. This suggests, that there is a limit to the geometry accuracy these approach.
If we look at the domains resulting from our reconstruction approaches on mesh levels four and five in \Cref{fig_stoke_domain_tip}, we see a visible discrepancy between the exact and reconstructed domains for all approaches. In particular, we see that the level-set approaches are qualitatively not as good as the spline approximation near the tip of the crack, which is consistent with the convergence results.

\subsection{Sneddon's Test Coupled to a Stationary Fluid-Structure Interaction Problem.}
\label{sec_test_FSI}
We now consider the stationary fluid-structure interaction 
problem \eqref{eqn.ale-fsi}, i.e., Formulation \ref{fsi:ale:stationary}, 
with the fluid domain $\FLref=\CR$ and $\SOref =\Omega\setminus\CR$, as reconstructed using the spline geometry approach. The force acting on the fluid is given by $\fbref\in L^2(\FLref)$ with the specific form 
\begin{equation}\label{eqn.fsi-forcing}
  \fbref = 
    \begin{pmatrix}
      \fbref_1 \\
      \fbref_2
    \end{pmatrix}
  = 
    \begin{pmatrix}
      0\\
      c_1 \exp(-c_2 \Vert\xbref- \xbref_0\Vert_2^2)
    \end{pmatrix}
\end{equation}
with the point source $\xbref_0 = (2.05, 2.01053) \in \FLref$, $\xbref\in \Oref$
and constants $c_1=10^{-4}$ and $c_2=10^3$. We note that $\xbref_0$ is located 
in a non-symmetrically in $\hat\FLref$.

As a quantity of interest, we look at the value of the deformation $\hat\ub$ at the point $\xb_{ref}=(2.1, 2.015795)$ which is close to the interface but inside the solid domain. As a reference value, we compute the problem on a domain constructed using the analytic values from Sneddon's test using a series of highly resolved meshes and high-order finite elements. We have found $\hat\ub(\xb_{ref}) = (-3.555\times10^{-11},1.303\times10^{-9})$ to be accurate to four significant figures.
 
\subsubsection{Results}

\begin{figure}
  \centering
  \begin{minipage}[b]{226pt}
    \centering
    \includegraphics[height=150pt]{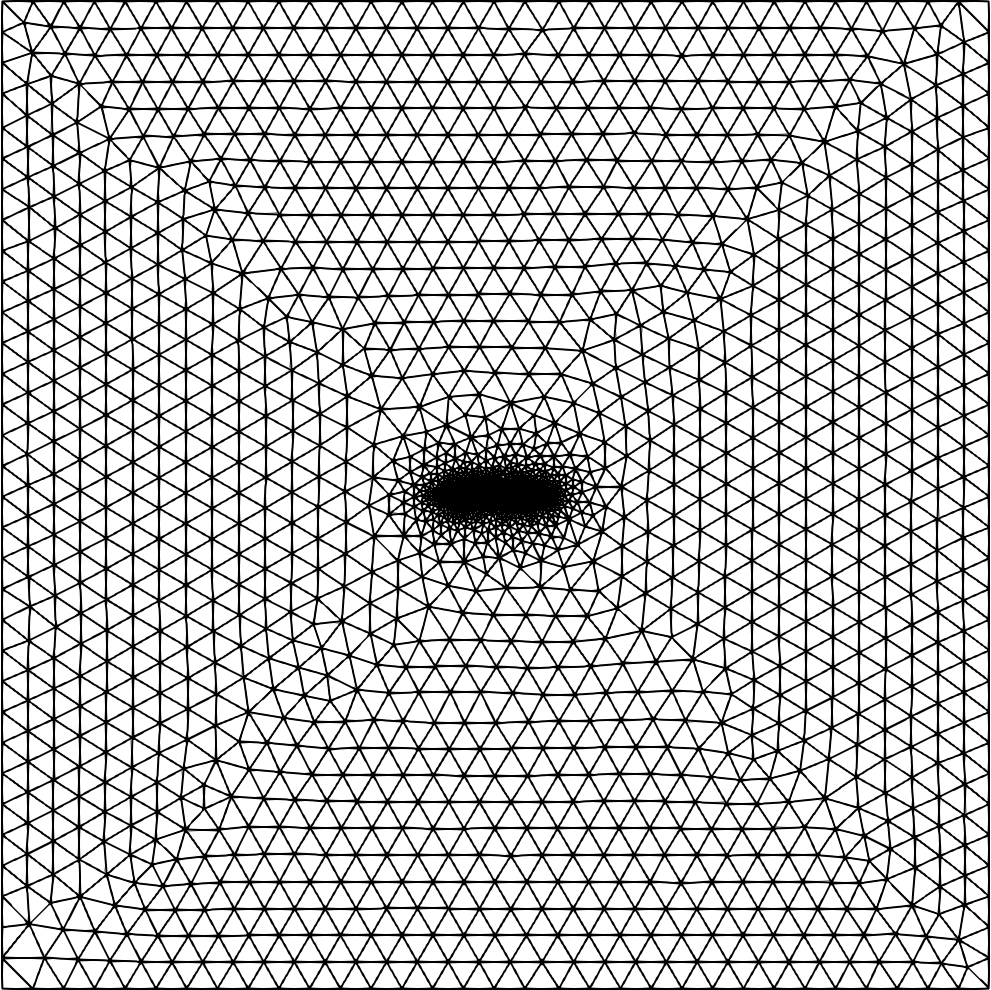}
  \end{minipage}
  \begin{minipage}[b]{226pt}
    \centering
    \includegraphics[height=150pt]{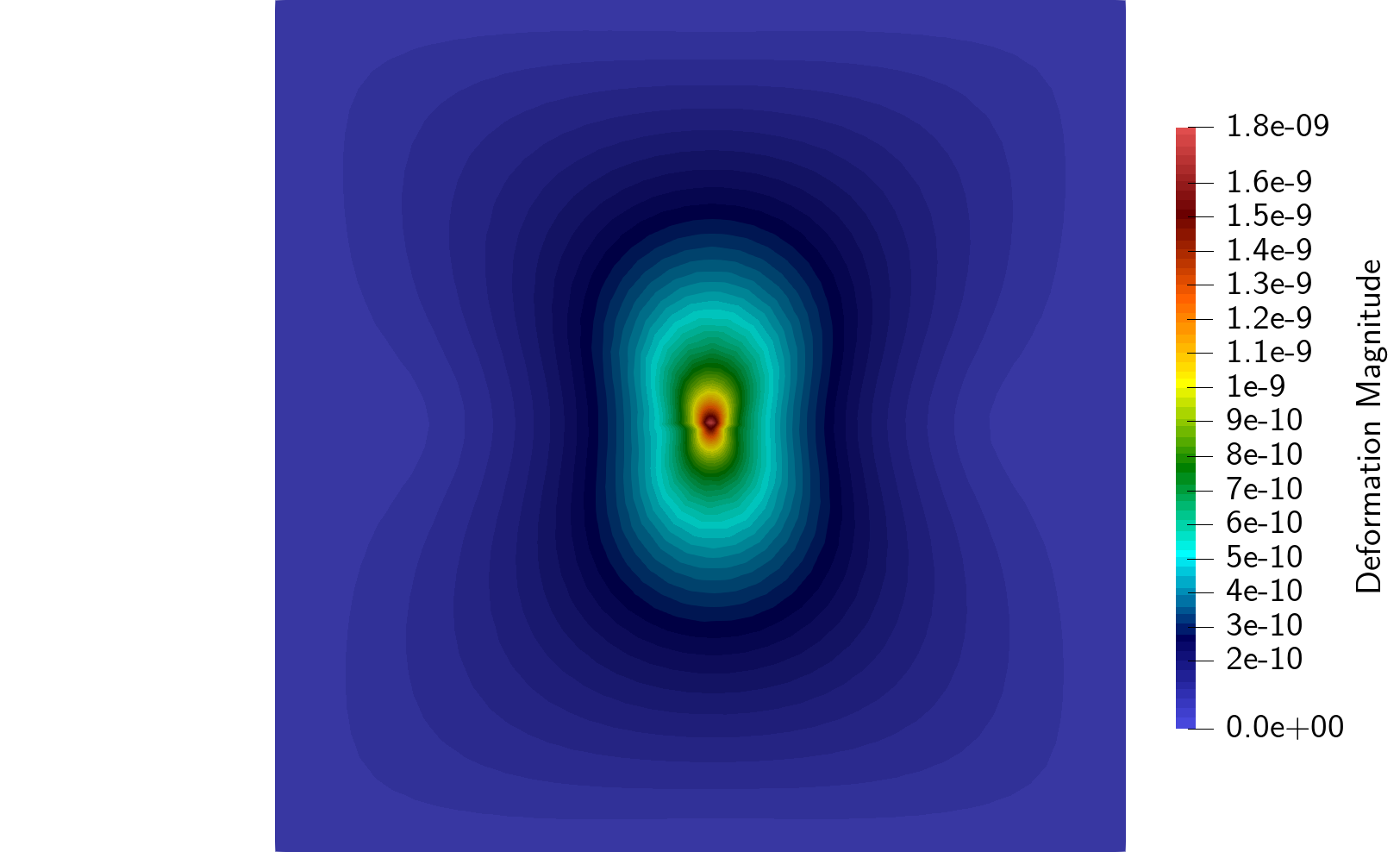}
  \end{minipage}
  
  \vspace{8pt}

  \begin{minipage}[b]{226pt}
    \centering
    \includegraphics[width=223.5pt]{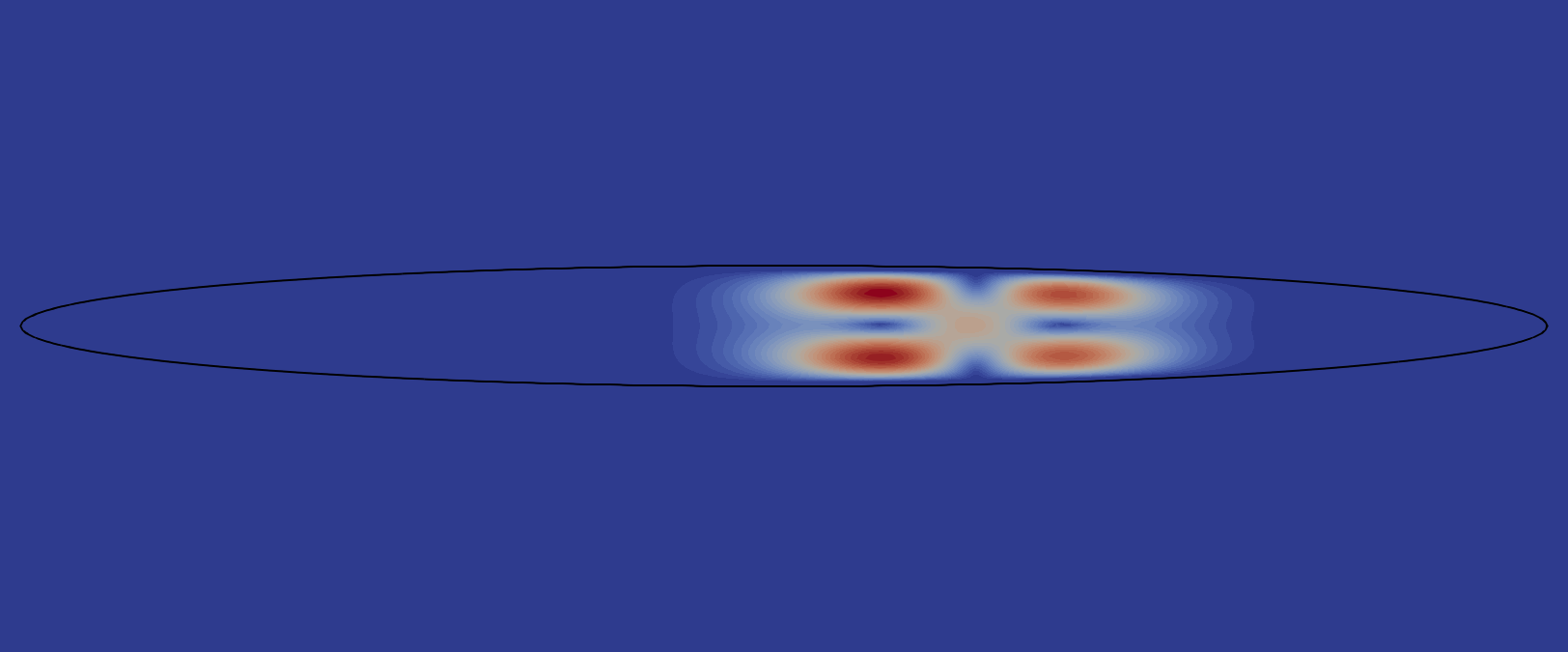}
    \includegraphics[width=223.5pt]{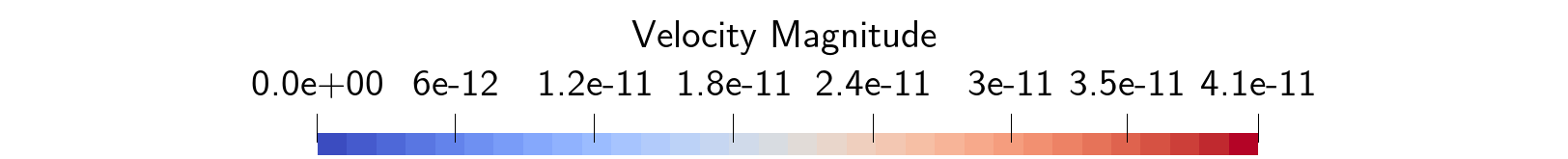}
  \end{minipage}
  \begin{minipage}[b]{226pt}
    \centering
    \includegraphics[width=223.5pt]{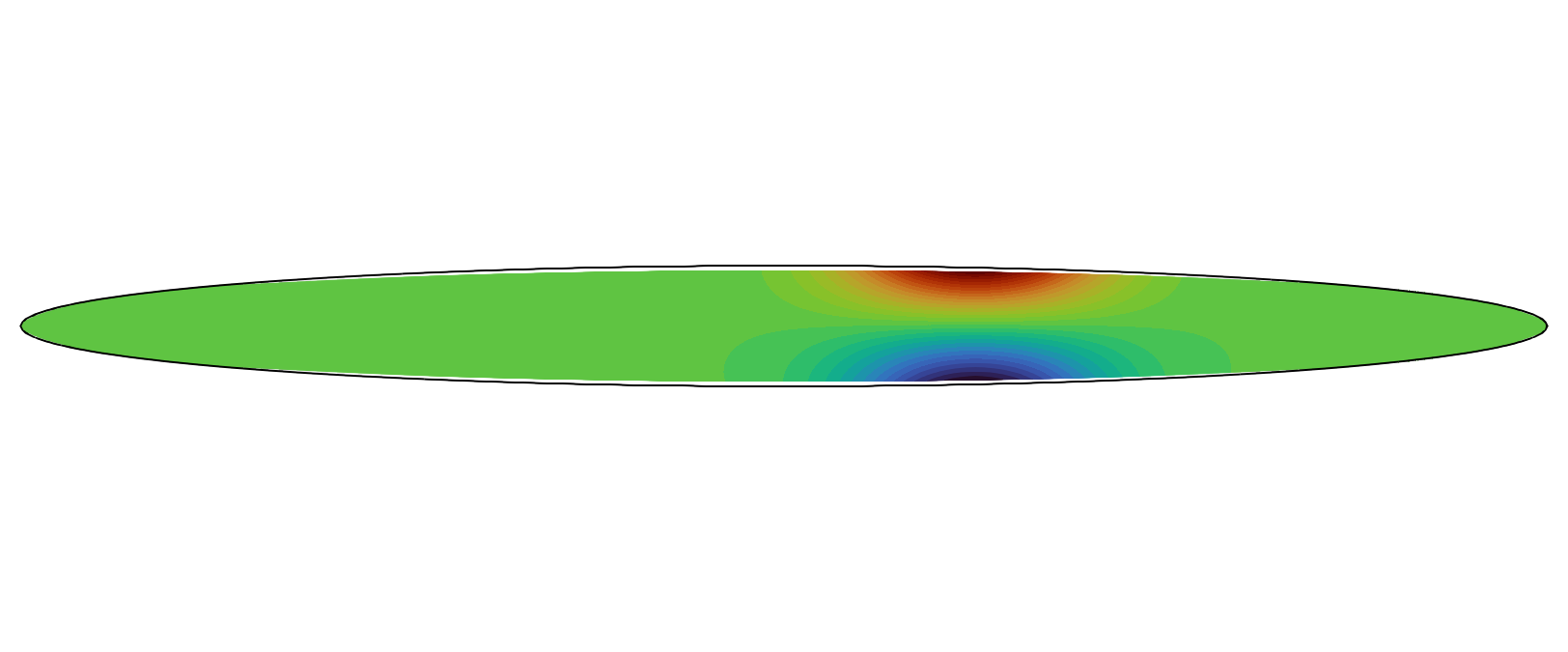}
    \includegraphics[width=223.5pt]{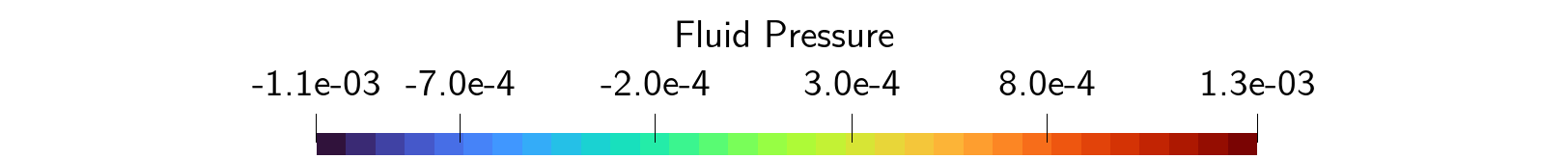}
  \end{minipage}
  \caption{Example \ref{sec_test_FSI}. Top left: Computational mesh of the reconstructed domain. Top Right: Deformation field. Bottom left: Velocity field zoomed to the open crack. Bottom Left: Pressure field inside the reconstructed crack.}
  \label{fig.fsi-results}
\end{figure}

We consider the problem over five meshes constructed from the crack opening displacements. The results can be seen in \Cref{tab.results.fsi}. Furthermore, the mesh, resulting displacement field, velocity and pressure can be seen in \Cref{fig.fsi-results}. We see that the point evaluation of the displacement field converges towards the expected value. However, the rate of convergence appears to be slow. We attribute this to the fact that the values are very small in absolute terms and, therefore, very sensitive to inaccuracies of the crack geometry. Furthermore, the analytical shape of the ellipse is derived under the assumption of an unbounded domain so that the exact size of the ellipse cannot be realised asymptotically on a bounded computational domain.

\begin{table}
  \centering
  \begin{tabular}{rll}
    \toprule
    Mesh level & $\hat\ub_1(2.1, 2.015795)$ & $\hat\ub_2(2.1, 2.015795)$\\
    \midrule
    0 & $-2.08958\times10^{-11}$ & $1.11351\times 10^{-9}$\\
    1 & $-2.88528\times10^{-11}$ & $1.21810\times 10^{-9}$\\
    2 & $-3.00343\times10^{-11}$ & $1.22641\times 10^{-9}$\\
    3 & $-3.15208\times10^{-11}$ & $1.24845\times 10^{-9}$\\
    4 & $-3.42140\times10^{-11}$ & $1.28361\times 10^{-9}$\\
    \cmidrule{1-3}
    Ref. & $-3.555\times10^{-11}$ & $1.303\times10^{-9}$\\
    \bottomrule
  \end{tabular}
  \caption{Example \ref{sec_test_FSI}. Results for the stationary fluid-structure interaction problem in the geometry reconstructed from the phase-field computation for Sneddons problem.}
  \label{tab.results.fsi}
\end{table}

\subsection{Two Orthogonal Cracks Coupled to a Stationary Fluid-Structure Interaction Problem}
\label{sec.text_fsi_T}

As a final test and second example for the full geometry reconstruction and fluid-structure interaction algorithm, we consider a more involved setting containing two orthogonal cracks. While the consideration of multiple cracks is not challenging for the phase-field computation itself; see, for example, the literature cited in the introduction and the numerous references therein. It is, however, more challenging with regard to the COD than the previous section, since multiple coordinate directions are involved. Consequently, the orthogonal crack will appear as a jump in the crack opening displacement of the first crack. As a result, the geometry reconstruction and meshing of the resulting geometry is also more involved, and this test serves as proof of concept for multiple fractures and more complex crack geometries.

\begin{figure}
  \centering
  \includegraphics{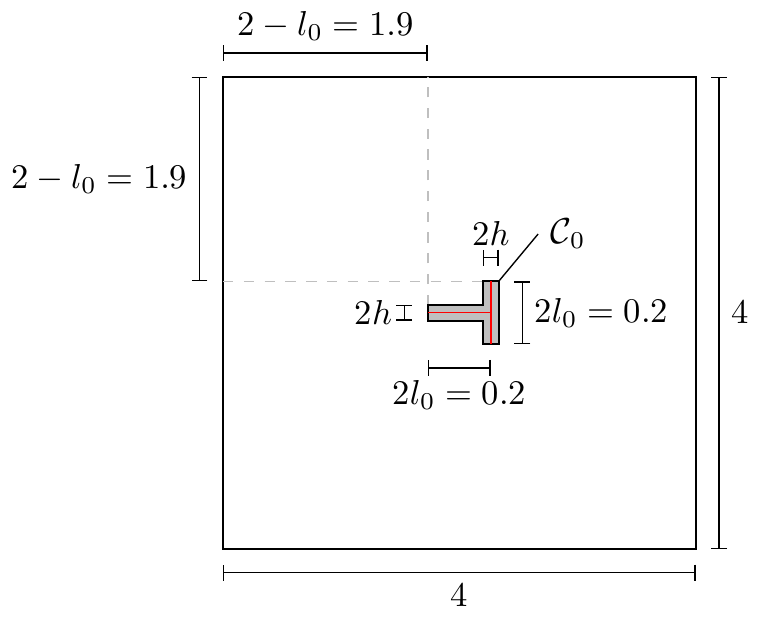}
  \caption{Example \ref{sec.text_fsi_T}. Schematics of the Phase-Field set-up.}
  \label{fig.schematics_fsi_T}
\end{figure}

We again consider the domain $\Omega=(0,4)^2$. The boundary conditions are again given by \eqref{eqn.senddon_bc1} and \eqref{eqn.senddon_bc2}. The initial phase field is a rotated ``T'' shape given by \eqref{eqn.sneddon_initial_pf} with
\begin{equation*}
  \CR_0 = (1.9, 2.1)\times(2-h, 2+h)\cup(2.1 - h, 2.1 + h)\times(1.9, 2.1),
\end{equation*}
i.e., two cracks of length $2l_0=0.2$. A sketch of this initial geometry can be seen in \Cref{fig.schematics_fsi_T}. The mechanical parameters are again Young's modulus and Poisson's ratio which we set in this case to be $E_s=5\times10^{4}$ and $\nu_s=0.35$. The applied pressure is $p=10^{4}$ and the critical energy release rate is chosen as $G_c=500$.

For the FSI problem, we consider the forcing term given in \eqref{eqn.fsi-forcing}, and $\xb_0=(2.098, 2.002), c_1=0.0001$ and $c_2 = 5000$, respectively. Consequently, the forcing is applied most strongly just off the intersection of the centrelines of the two cracks. As a specific quantity of interest, we will consider the solid deformation in the point $\xb=(2.05, 2.025)$.

\subsubsection{Results}
\begin{figure}
  \centering
  \includegraphics[width=200pt]{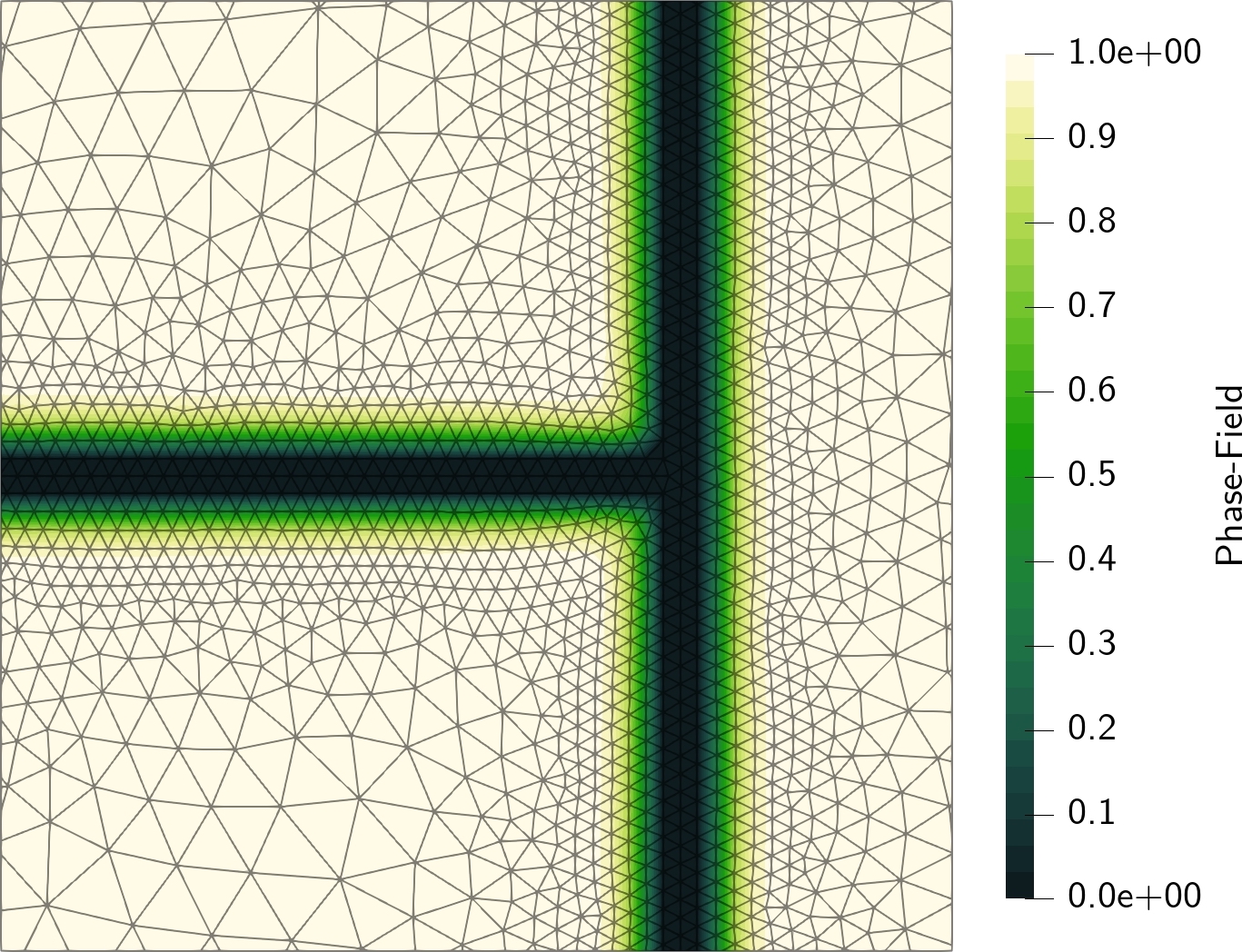}
  \hspace*{.5cm}
  \includegraphics[width=200pt]{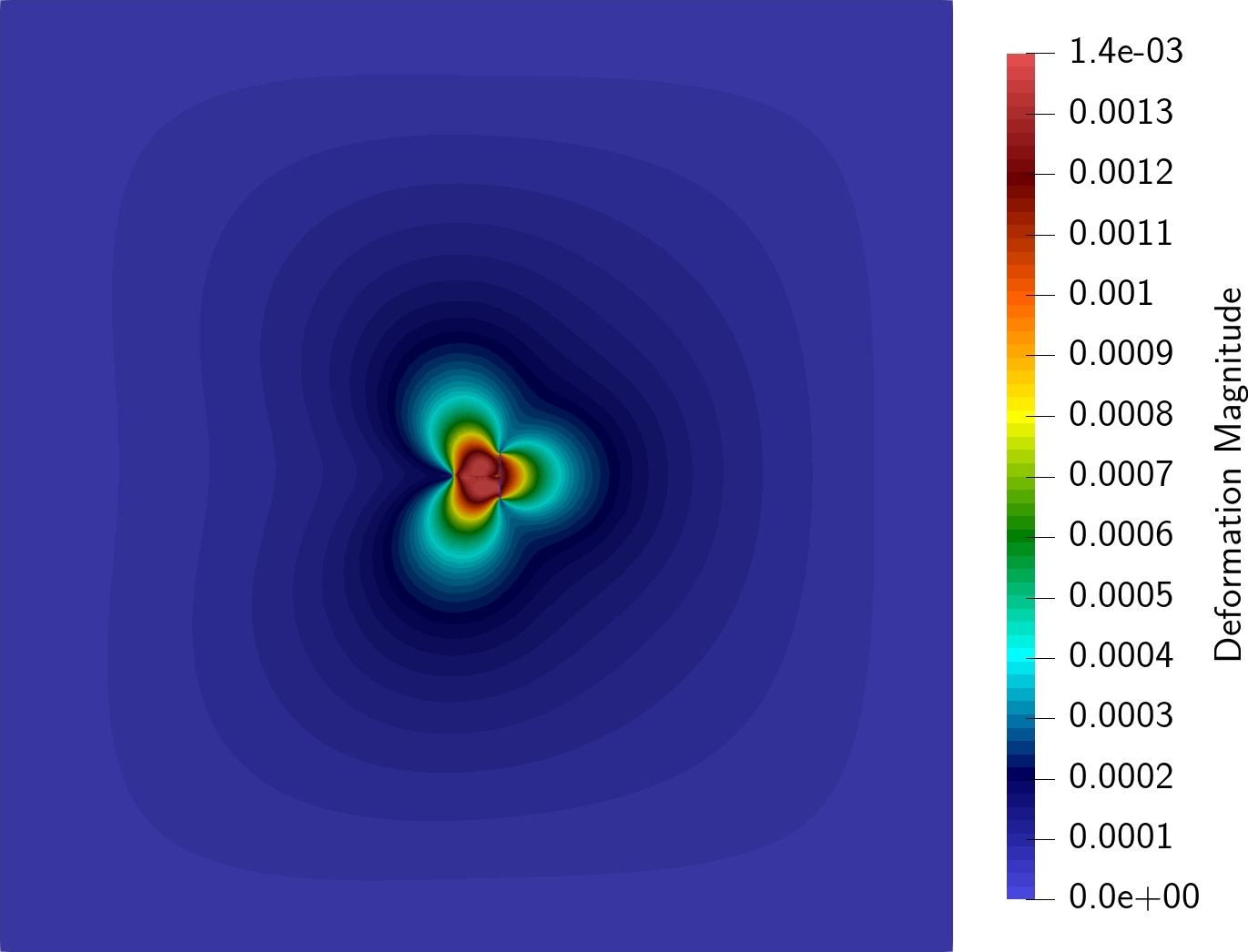}
  \caption{Example~\ref{sec.text_fsi_T}. Left: Phase-field solution and computational mesh zoomed-in to $(2.075,2.11)\times(1.9825,2.0175)$. Right: Phase-field deformation solution on the entire domain. Computed on mesh level 4.}
  \label{fig.phase-field_fsi_T}
\end{figure}

\begin{figure}
  \centering
  \includegraphics[width=200pt]{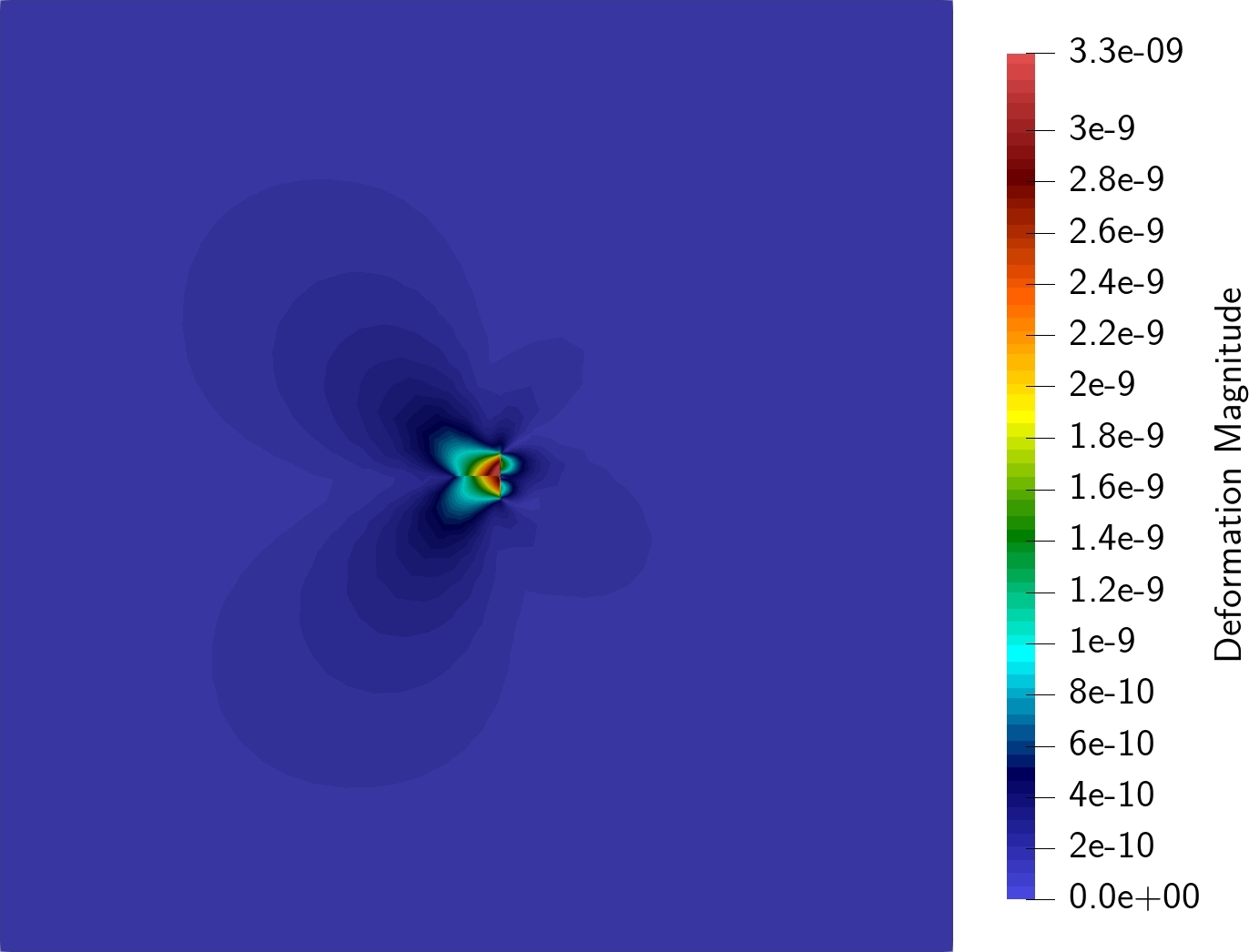}
  \hspace{.5cm}
  \includegraphics[width=200pt]{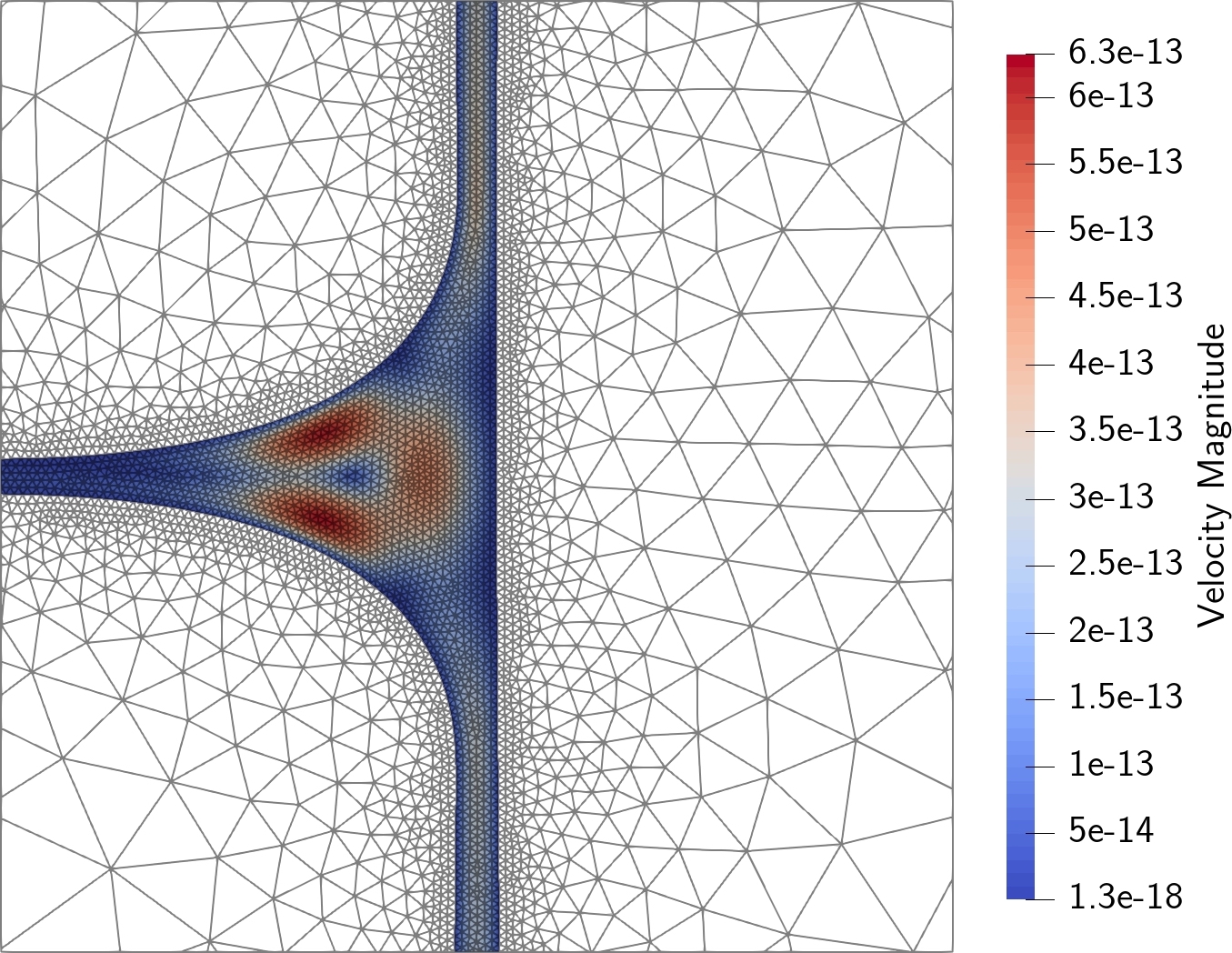}
  \caption{Example~\ref{sec.text_fsi_T}. Left: Solid deformation in the entire domain. Right: Fluid velocity solution zoomed to the region $(2.08, 2.12)\times(1.98, 2.02)$, together with the mesh of the reconstructed domain. Computed on mesh level 4.}
  \label{fig.deform_vel_fsi_T}
\end{figure}

We consider the initial mesh size $h_0 = 0.01$ and 5 levels of mesh refinement. The resulting phase-field deformation and a zoom-in of the phase-field solution on mesh level four can be seen in \Cref{fig.phase-field_fsi_T}. The resulting point evaluation of the FSI deformation is presented in \Cref{tab.restults_fsi_T}. The FSI deformation on mesh level four is shown on the left of \Cref{fig.deform_vel_fsi_T}, and a zoom-in of the velocity solution and the FSI mesh around the crack intersection can be seen on the right of \Cref{fig.deform_vel_fsi_T}. The FSI deformation results are consistent and converge towards approximately $(-2\times 10^{-9}, -3.9\times 10^{-10})$, showing that the domain reconstruction works. Looking at the FSI deformation, we see that this is largest in the region to the top-left of the crack intersection. This is to be expected since the forcing term's peak in the fluid domain has been positioned to the top left of the crack intersection. 

\begin{table}
  \centering
  \begin{tabular}{rll}
    \toprule
    Mesh level & $\hat\ub_1(2.098, 2.0002)$ & $\hat\ub_2(2.098, 2.0002)$\\
    \midrule
    0 & $-2.39501\times 10^{-9}$ & $\phantom{-}2.38784\times 10^{-10}$\\
    1 & $-2.88056\times 10^{-9}$ & $-1.04270\times 10^{-10}$\\
    2 & $-3.03558\times 10^{-9}$ & $-6.25539\times 10^{-10}$\\
    3 & $-2.78315\times 10^{-9}$ & $-6.23526\times 10^{-10}$\\
    4 & $-2.60130\times 10^{-9}$ & $-4.75054\times 10^{-10}$\\
    5 & $-2.33500\times 10^{-9}$ & $-3.96828\times 10^{-10}$\\
    \bottomrule
  \end{tabular}
  \caption{Example~\ref{sec.text_fsi_T}. Results for the stationary fluid-structure interaction problem in the geometry reconstructed from the phase-field of two orthogonal cracks.}
  \label{tab.restults_fsi_T}
\end{table}

\section{Conclusions}
\label{sec_conclusions}

In this work, we proposed a high-accuracy framework for multi-domain multi-physics phase-field fracture. The main drawback of phase-field fracture formulations is a smeared transition zone which results in a significant loss of accuracy when interface conditions are described. To this end, we proposed that phase-field fracture serves as a predictor for the crack path, followed by a reconstruction of the fracture geometry. With this procedure at hand, we arrived from an interface-capturing method
at an interface-tracking scheme in which interface conditions can be described with high accuracy. This is realised in one algorithm and one software, namely \texttt{NGSolve}, and therefore a promising procedure. As an application, we considered fluid-filled fractures. Therein, a Stokes flow was described inside the fracture, which then was coupled in a one-way procedure to the surrounding elastic medium. To demonstrate the robustness and accuracy of our framework, we first investigated Sneddon's benchmark problem, which is a well-accepted example. Here, we obtained excellent findings. Next, we considered Stokes flow in the fracture only, again with convincing findings. 
In our third numerical example, we coupled the fracture Stokes flow to the surrounding elastic medium. This resulted in a fluid-structure interaction problem, which was treated with the interface-tracking arbitrary Lagrangian-Eulerian approach. 
For our fourth test, we considered a more involved example for the full geometry reconstruction and stationary Stokes coupled to a elastic medium algorithm. This test consisted of two orthogonal cracks, making the geometry reconstruction more challenging as it is no longer aligned to one single coordinate axis. Moreover, this last example shows that our procedure is not restricted to single fractures.
Our results suggest a promising procedure for such complex problem statements when a high accuracy at the interface is indispensable.

We note that the next future extension is to create an iterative loop taking the Stokes pressure and computing again the pressurised phase-field problem, followed by a new reconstruction of the interface, and then computing the fluid-structure interaction problem again.

\section*{Data Availability Statement}
The source code and the data generated with it is publicly available on github under \url{https://github.com/hvonwah/stationary_phase_field_stokes_fsi} and archived on zenodo \cite{vWW22_zenodo} under the GNU General Public License v3.0.

\section*{Acknowledgments}
HvW acknowledges support through the Austrian Science Fund (FWF) project F65.

\printbibliography

\end{document}